\documentclass[a4paper, english, 11pt]{article}

\usepackage{setspace}
\usepackage[T1]{fontenc}
\usepackage{hyperref}
\usepackage[utf8]{inputenc}
\usepackage{dcolumn}
\usepackage{color}
\usepackage[x11names]{xcolor}
\usepackage[thmmarks]{ntheorem}
\usepackage{ragged2e}
\usepackage{pdfpages}
\usepackage{amsmath}
\usepackage[french, english]{babel}
\usepackage{multirow}
\usepackage{array}
\usepackage{booktabs}
\usepackage{wrapfig}
\usepackage{multicol}
\usepackage{algorithm}
\usepackage{lscape}
\usepackage{pifont}
\usepackage{graphicx}
\usepackage{cite}
\usepackage{calc}

\usepackage[all]{xy}
\usepackage{framed}
\usepackage{tikz}
\def\dar[#1#2]{\ar@<2pt>[#1]\ar@<-2pt>[#2]}

\usepackage{amsfonts,amssymb,amsmath}
\usepackage[top=3cm,bottom=3cm,left=3cm,right=3cm]{geometry}
\usepackage{fancyhdr}
\usepackage{array}

\newcommand{\kk}[0]{\textbf{k}}
\newcommand{\Mod}[0]{\text{Mod}}

\newcommand{\Pe}{\text{Pers}}

\newcommand{\R}[0]{\mathbb{R}}
\newcommand{\N}[0]{\mathbb{N}}

\newcommand{\Z}[0]{\mathbb{Z}}
\newcommand{\C}[0]{\mathcal{C}}

\newcommand{\B}[0]{\mathcal{B}}

\newcommand{\V}[0]{\mathbb{V}}
\newcommand{\D}[0]{\text{D}}
\newcommand{\Obj}[0]{\text{Obj}}
\newcommand{\Ho}[0]{\text{H}}

\newcommand{\Ouv}[0]{\mathrm{Open}}
\newcommand{\op}[0]{\text{op}}
\newcommand{\Hom}[0]{\text{Hom}}
\newcommand{\Pers}[0]{\text{Pers}}
\newcommand{\Rr}[0]{\text{R}}
\newcommand{\MV}{\text{MV}}
\newcommand{\s}{\textbf{sf}}

\newcounter{nfigure} 
\makeatletter
 
\newskip\@bigflushglue \@bigflushglue = -100pt plus 1fil

\def\bigcentering{\let\\\@centercr\rightskip\@bigflushglue%
\leftskip\@bigflushglue
\parindent\z@\parfillskip\z@skip}

\theoremheaderfont{\scshape\bfseries}
\theorembodyfont{\normalfont}
\newtheorem{prop}{Proposition}[section] 

{ \theoremsymbol{\flushright{$\square$}} 
 \newtheorem*{pf}{Proof.} }
 { \theoremsymbol{\flushright{$\square$}} 
 \newtheorem*{pfofLemma}{Proof of the lemma.} }
   
 \newtheorem{ex}[prop]{Example} 
 \newtheorem{cor}[prop]{Corollary}
 \newtheorem{remark}[prop]{Remark}
 \newtheorem{lem}[prop]{Lemma}
 \newtheorem{defi}[prop]{Definition}
 \newtheorem{thm}[prop]{Theorem}
\author{Nicolas Berkouk, Grégory Ginot, Steve Oudot}

\begin{document}
\selectlanguage{english}

\title{Level-sets persistence and sheaf theory}
\maketitle
\begin{abstract}
   In this paper we provide an explicit connection between  level-sets persistence and derived sheaf theory over the real line. In particular we construct a functor from 2-parameter persistence modules to sheaves over $\R$, as well as a functor in the other direction. We also observe that the 2-parameter persistence modules arising from the level sets of Morse functions carry extra structure that we call a Mayer-Vietoris system. We prove classification, barcode decomposition, and stability theorems for these Mayer-Vietoris systems, and we show that the aforementioned functors establish a pseudo-isometric equivalence of categories between derived constructible sheaves with the convolution or (derived) bottleneck distance and the interleaving distance of strictly pointwise finite-dimensional Mayer-Vietoris systems. Ultimately, our results provide a functorial equivalence between level-sets persistence and derived pushforward for continuous real-valued functions.
\end{abstract}

\tableofcontents

\section{Introduction}

Persistent homology is a powerful and versatile tool of applied algebraic topology that has found applications in a variety of areas of the Sciences, in particular those connected with data Science. 
Roughly speaking, the aim of persistent homology is to define algebraic invariants for filtered topological spaces that are both robust and computer-friendly. Typically, given a continuous function $f: X\to \R$, one  considers the homology of the \emph{sublevel sets} $H_i(f^{-1}((-\infty, t)))$, for a real parameter $t$, and of the various inclusion maps induced between them as $t$ grows. This data is called the \emph{sublevel-sets persistence module} associated to $f$. Under some reasonable finiteness assumptions, this data can be compactly encoded with a set of intervals---called a {\em barcode}---that describe when the homology generators appear, live, then die in the family of sublevel sets.
Barcodes are easy to compute and to handle on a computer. They serve as descriptors for data in applications, where they can be compared  using a matching distance called the {\em bottleneck metric}. This metric is actually  equivalent to the so-called  {\em interleaving distance} between persistence modules, which is both {\em stable} (i.e. robust to perturbations of the input function~$f$) and {\em universal} (i.e. the most sensitive among stable distances)~\cite{Lesn12}. These properties constitute  the cornerstone of persistence theory and they play a key role in its applications. 

\medskip 

In the recent years, several generalizations of persistent homology have been proposed, whose goal is to extract richer information on  the topology of a continuous function $f:X\to \R$. 
In this work we focus on two of the most prominent ones, namely {\em level-sets persistence} and {\em two-parameter persistence}, and we provide an explicit connection to derived (co)sheaf theory.

Level-sets persistence studies the homology groups of preimages
$H_i(f^{-1}(]s,t[)$, where $]s,t[$ is the French notation\footnote{We
      adopt this notation for the sake of clarity, to avoid potential
      confusions with the point~$(s,t)\in\R^2$.} for the open
    interval~$(s,t)$. The collection of these groups for $s<t\in\R$,
    together with the collection of morphisms induced by inclusions of
    smaller intervals into larger intervals, form a two-parameter
    persistence module indexed over the upper half-plane $\Delta^+ =
    \{(x,y) \mid x+y>0\}$ via the identification of each
    interval~$]s,t[$ with the point $(-s,t)$. This module is called
    the {\em level-sets persistence module} of $f$ and denoted
    by~$M^f$ (see~\ref{Ex:MVfromFct}).  Note that the plane $\R^2$ is
    equipped with the partial product order, noted~$\leq$. Henceforth
    we will write $\Pers(\R^2)$ for the category of persistence
    modules indexed over~$(\R^2, \leq)$, and $\Pers(\Delta^+)$ for its
    counterpart over $(\Delta^+, \leq)$.  Unfortunately, though very
    natural, the theory of general 2-parameter persistence modules is
    significantly more complicated than 1-parameter persistence. For
    instance, there is no analogue of barcodes, due to the
    poset $(\R^2,\leq)$ being a wild-type quiver with arbitrarily
    complicated indecomposables.
    Nevertheless, one can still
    define an interleaving distance in this context,  satisfying the same  stability and universality properties as in 1-parameter persistence~\cite{Lesn12}.

\medskip 

Another very promising direction of investigation is given by merging (co)sheaves theory with persistence and computer-friendly techniques. It was pioneered by the work of Curry~\cite{Curr14}, and a general framework was developed by Kashiwara and Schapira~\cite{Kash18, KS18}. In order to benefit fully from the algebraic topology of sheaves, it is necessary to work with the derived (sometimes called homotopy) category. Kashiwara-Schapira have equipped the derived category of sheaves with a distance, called the \emph{convolution distance}, which is a derived analogue of the interleaving distance, see~\cite{Kash18, BerkoukPetit}. Furthermore, there is a natural notion of barcode and a decomposition theorem for \emph{constructible} sheaves over $\R$. To a function $f:X\to \R$, one can associate a canonical sheaf over $\R$, namely the derived pushforward $\Rr f_* \kk_X$, which is a sheaf analogue of the level-sets persistence homology introduced earlier (here $\kk$ is our ground field). 
The persistence theory for sheaves over $\R$, following Kashiwara and Schapira's program, was investigated in depth in~\cite{Berk18}. In particular,  a \emph{derived bottleneck distance} for constructible sheaves was developed and proven to be \emph{isometric to the convolution distance}.  

\medskip 

The main motivation of this paper is to relate \emph{precisely} the above developments, namely: persistence modules over $\Delta^+:=\{(x,y),\, x+y>0\}$ on the one hand; sheaves over $\R$ on the other hand. Specifically, given a function $f:X\to\R$, we are interested in connecting  the level-set persistence module $M^f$ with the derived pushforward $\Rr f_* \kk_X$. 
In order to do so, we will construct a functor~$\Xi$ from 2-parameter persistence modules to sheaves over $\R$ (see section~\ref{SS:levelsettopresheaves}). 
Note that this functor is not an equivalence of categories, and that it is not isometric nor reasonably Lipschitz either. Indeed, there can be no equivalences or almost equivalences between these two categories, since the general category of 2-parameter persistence modules is wild representation type as we have mentioned already, and since its objects do not, in general, satisfy any of the local-to-global properties of sheaves.  

\smallskip

As mentioned above, of particular interest to us are the level-sets persistence modules $M^f$ arising from continuous functions $f:X\to \R$, which actually have  \emph{more structure} than general persistence modules over $\Delta^+$. 
\begin{center}
\emph{Our idea is thus to consider a variant of the category of 2-parameter persistence modules taking into account the extra structure and properties carried by level-sets persistence modules.}\end{center} This follows the fundamental credo of algebraic topology that extra structure on homology gives refined homotopy and geometric information. A general idea here is that to get a better-behaved category of 2-parameter persistence modules, it is key to consider and restrict to those objects having the extra structure and properties coming from  data arising in practical applications.

\smallskip

Let us now explain where this extra structure comes  from: the  various homology groups of a topological space obtained as the union of two open subsets  are connected through the well-known Mayer-Vietoris long exact sequence. This sequence involves the homology groups of the union, the sum of the homology groups of the two open  subsets, and the homology of their intersection. 
We axiomatize this data to define  a structure we call \emph{Mayer-Vietoris (MV) persistence systems over $\Delta^+$}, whose category is denoted by $\text{M-V}(\R)$. 
A MV-system is a graded persistence module $(S_i)_{i\in \Z}$ over $\Delta^+$,   together with connecting morphisms $\delta^s_i: S_{i+1}[s]\to S_i$ for all vectors $s\in\Delta^+$ and grades $i\in\Z$, giving rise to the following  \emph{exact} sequences  (see Definition~\ref{D:MVSystem}):
\[\xymatrix{S_{i+1}[s] \ar[r]^-{\delta_{i+1}^s} & S_i \ar[r] & S_i[s_x] \oplus S_i[s_y] \ar[r] & S_i[s] \ar[r]^-{\delta_{i}^s} & S_{i-1}}\]  

and satisfying some appropriate compatibility conditions. 
These sequences encode  the interactions between the various homology groups at various points of $\Delta^+$, and they carry both a derived and local-to-global information---in some sense that will be made precise in the paper. 

A key property that we leverage in our analysis, is that the category of Mayer-Vietoris systems is rather well behaved. In particular, we prove a \emph{structure theorem} for Mayer-Vietoris persistence systems under standard pointwise finite dimensionality assumptions (see  Theorem~\ref{thm:decom_MV}). According to this result, there are \emph{four different types} of \emph{indecomposable} Mayer-Vietoris systems, which all have pointwise dimension at most~1 and are therefore characterized by their supports. The supports can be either vertical or horizontal bands, or else birth or death blocks (see Definition~\ref{D:blocksmodulesforMV} and Lemma~\ref{L:Blockassociatedtointerval}). Degree-wise, these indecomposables behave like the so-called \emph{block modules} from level-sets persistence and middle-exact bipersistence theories~\cite{BotCra18,BoLe16,CdSKM19,CdSM09,CO17}. For this reason, in the following we abuse terms and also call our indecomposables \emph{block MV-systems}. 
Our structure theorem (Theorem~\ref{thm:decom_MV}) takes the following form:
\begin{thm}
{A, bounded below, pointwise finite-dimensional (pfd) Mayer-Vietoris system has a unique decomposition as a direct sum of block MV-systems.}
\end{thm}
This result follows non-trivially from the decomposition theorem for middle-exact bipersistence modules~\cite{BotCra18,CO17}. It provides a \emph{barcode for Mayer-Vietoris systems}, made of the blocks involved in their decomposition. Furthermore, we have a canonical \emph{interleaving distance for Mayer-Vietoris systems} since they form a category.

\smallskip 

The aforementioned functor~$\Xi$ from 2-parameter persistence  modules to sheaves lifts as a (contravariant) functor $\overline{(-)}^{MV}$ from Mayer-Vietoris systems to the derived category $\D(\kk_\R)$ of sheaves on $\R$, which is essentially the sheafification of the duality functor. We \emph{construct a pointwise section} of this functor, i.e. a functor $\Psi$ from sheaves to Mayer-Vietoris systems such that the composition with $\overline{(-)}^{MV}$ gives the identity pointwise on every sheaf $F\in \D^b_{\R c}(\kk_\R)$ (see Corollary~\ref{C:MVcircPsiisId}): $$\left ((\overline{~\cdot~ })^{\text{MV}} \circ \Psi \right ) (F) \simeq F.$$
Roughly speaking, this functor $\Psi$ is defined as the dual of the derived global sections of sheaves  (see Definition~\ref{def:Psi}). Both functors restrict to the subcategories of pointwise (resp. strictly pointwise see Definition~\ref{D:spfd}) finite-dimensional Mayer-Vietoris systems on one side, and of constructible sheaves on $\R$ on the other side. 
Under standard  pointwise finiteness conditions, we are able to prove that these two functors establish a pseudo-isometric equivalence between these categories. 
More precisely, our second main theorem (see Theorem~\ref{T:MainIsometry} and Corollary~\ref{C:MVdisO}) states as follows: 
\begin{thm}\label{T:MainIsoIntro} The functors $\overline{(-)}^{MV}$ and $\Psi$ form a pseudo-isometric equivalence of categories, meaning:
\begin{itemize}\item for all strictly pointwise finite-dimensional Mayer-Vietoris systems $M, N$,  one has equality 
$d_I(M, N) = d_C (\overline{M}^{MV}, \overline{N}^{MV}) = d_B(\mathcal{B}(\overline{M}^{MV}), \mathcal{B}(\overline{N}^{MV})) $ between the interleaving, convolution and derived bottleneck distances; 
\item for all constructible sheaves $F, G \in \D^b_{\R c}(\kk_\R)$, one has 
$ d_B(\mathcal{B}(F), \mathcal{B}(G)) = d_C(F,G) =d_I(\Psi(F), \Psi(G))$;
\item  $\overline{M}^{MV} = \overline{N}^{MV} $ if and only if $d_I({M}, {N})=0$.\end{itemize}
\end{thm}
In particular, the derived distances can be computed using the 2-paramater interleaving distance for M-V systems.


To prove this theorem, in Sections~\ref{SS:MVfunctor} and \ref{SS:Psi} we explicitly compute the action of the sheafification of Mayer-Vietoris systems functor $\overline{(-)}^{MV}$ and of its section $\Psi$ on shifts and convolution for the building block modules of each theories. These  are  computations of independent interest.

\smallskip

Finally, we prove that the functors $\Rr f_* \kk_X$ and $M^f$ are equivalent to each other under these transformations, i.e. $\overline{M^f}^{MV} \cong (\Rr f_* \kk_X)$ (Proposition~\ref{P:MVmapstoRf}), hence isometric according to theorem~\ref{T:MainIsoIntro}, thus establishing the sought-for correspondence between level-sets persistence and derived pushforward for continuous real-valued functions. 

\smallskip

Theorem~\ref{T:MainDiag} below summarizes our main results connecting level-sets persistence, Mayer-Vietoris systems, and derived sheaves. The notation $\mathrm{Top}_{|\R}$ stands for the category of topological spaces over $\R$, whose objects are spaces $X$ together with a continuous map  $f:X\to\R$, and whose morphisms are commutative triangles $\xymatrix{X\ar[r]_{\phi} \ar@/^1pc/[rr]^{f} & Y\ar[r]_{g} & \R } $. We let $ \mathrm{Top}^{c}_{|\R} $ denote the subcategory of those functions $f:X\to \R$ such that $\Rr f_* \kk_X$ is constructible,  $\D^b_{\R c}(\kk_\R)$ denote the bounded derived category of constructible sheaves on~$\R$, and we denote $\Rr(-)_*\kk_{(-)}$ the functor $f\mapsto \bigoplus \Rr^if_*(\kk_X)[-i]$.
\begin{thm}\label{T:MainDiag} The following diagram of categories and functors commutes:
\[\xymatrix{ \mathrm{Top}_{|\R} \ar[r]_{M^{(-)}} \ar@/^2pc/[rr]^{\Rr(-)_*\kk_{(-)}} &\text{M-V}(\R) \ar[r]_{\quad\overline{(-)}^{MV}} &  \D(\kk_\R)^{\op} \\ 
 \mathrm{Top}^{c}_{|\R} \ar[r]^{M^{(-)}} \ar@{^{(}->}[u] \ar@/_2pc/[rr]_{\Rr(-)_*\kk_{(-)}}&\text{M-V}(\R)^{\s} \ar@{^{(}->}[u]\ar[r]^{\quad\overline{(-)}^{MV}} & \D^b_{\R c}(\kk_\R)^{\op}. \ar@{^{(}->}[u] }  \] 
 Furthermore, $\overline{(-)}^{MV}$ and the vertical maps are isometries, while the other maps are $1$-Lipschitz.
\end{thm}
\begin{pf}
 The existence and commutativity of the diagram is the content of Lemma~\ref{L:MVofsfrestricts}, Propositions~\ref{P:MVfromFct} and \ref{P:MVmapstoRf}, and Theorem~\ref{T:MainIsometry}. The rest of the statement is given by Theorem~\ref{T:MainIsometry} and the stability theorems~\ref{P:StabilityfordMV},~\ref{P:StabforConv}.
\end{pf}
Theorem~\ref{T:MainDiag} relates precisely, and in fact essentially identifies, level-set persistence and constructible sheaves. Moreover, it does so in a functorial way.
In the final section of the paper (Section~\ref{SS:Projectionsfromcircle}) we give an example with full detail that illustrates this result. 

\subsection{Related and future work}

Parallel and independently to our work, Fluhr~\cite{Flu18} has consider similar problems, building a functor~$h$  from the category of derived constructible sheaves over $\R$ to the category of contravariant functors on some poset~$\mathcal{B}$. Once restricted to the subposet of all points/elements corresponding to bounded intervals, such contravariant functors carry roughly the same data as our Mayer-Vietoris systems.

In future work it will be interesting to use the isometric functors we have constructed to study the shriek functors associated to functions, namely $\Rr f_{!} \kk_X$. 

Furthermore, the pseudo-isometry theorem~\ref{T:MainIsoIntro} shall lift at the level of derived category of 2-paramater persistence modules to give an isometric equivalence of categories between constructible sheaves and  the quotient of a  full subcategory of  2-parameter persistent  complexes   where one identifies  those MV-systems that are at distance~$0$ of each other. 

\subsection{Notations}
\noindent  Here we detail our notations and conventions, for the reader's convenience:
\begin{itemize}
\item We fix a ground (commutative) field denoted $\kk$.
    \item For $s=(s_1,s_2)\in \R^2_{>0}$, we will use the notations $s_x:=(s_1,0)$ and $s_y := (0,s_2)$.
    \item Given a category $\C$, we denote $\C^{\op}$ its opposite category. 
    \item We will use the same notation for a poset $(S, \leq)$ and its associated  category whose objects are the elements of $S$ and whose set of morphims from  $s$ to $t$ consists of a single element if $s\leq t$ and is empty if not.
    \item A functor $M$ from a poset $(S, \leq)$ to vector spaces will be called \emph{pointwise finite dimensional}, \emph{pfd} for short, if for every $s\in S$, $M(s)$ is finite dimensional. 
    \item We denote $t\mapsto \vec{t}$ the functor $(\R ,\leq) \to (\Delta^+, \leq)$.
    \item For a death block $B$, see Definition~\ref{D:dualblocks}, we denote $B^\dag$ its dual  in $\R^2$, and vice-versa.
    \item The derived category associated to an abelian category $\mathcal{C}$ (such as complexes of sheaves or complexes of $\kk$-modules) is the localization of $\mathcal{C}$ with respect to quasi-isomorphisms, that is the category obtained by formally inverting quasi-isomorphisms. See e.g.~\cite{Kash90, Tohoku}.
\end{itemize}
\subsubsection{Notations for  intervals} 
\begin{itemize} \item We will use the \emph{french notation} $]a,b[$ for \emph{open intervals} $(a,b)$ in $\R$. The reason is to avoid confusion with points $(a,b)\in \R^2$ which will both be possible values of persistent or sheaf objects (and usually appear with similar letters).
    
\item For real numbers $a\leq b$, the notation $\langle a, b\rangle$ will mean an  interval whose boundary points are $a$ and $b$. We use this notation $\langle\, , \,  \rangle$  when we do not want to precise if the interval is open, compact, or half-open; in other words as a variable.
    \end{itemize}
\subsubsection{Notations  and conventions for shifts of graded and persistent objects}   
Standard and convenient notations for shifting the degree of a (possibly differential) graded object or for  shifted (or translated) persistent object are both given by  $[-]$ in the literature. We will have to use objects which are both differential graded and persistent, and we now explain how to avoid confusion about this notation in the paper. Note that we also have to deal with objects which are naturally homologically graded (for instance persistence modules) and cohomologically graded (sheaves).

For any (differential) cohomologically \emph{graded} object $C$, we will use the notation $C[i]$ for the (differential) graded object $C[i]^n:= C^{i+n}$ where $i\in \Z$. The letter $i$ can be replaced by $j$, $k$,  $\ell$, $m$ or $n$ in the paper, and the notation with one of these letters always means such a grading shift. These letters can also show up in subscripts. 

\smallskip

Similarly, for a (differential) \emph{homologically graded} object $M_*$ we will use the notation $M_*[i]$ for the graded object $(M_*[i])_n= M_{i+n}$, following for instance the conventions of~\cite[Section 12.13]{stacks-project}. We warn the reader that there is also  an opposite convention in the literature (which is the topological convention for suspension). The main advantage of this choice of convention in this paper is that the duality functor commutes with the shift in grading (instead of changing to its opposite):
$$\Hom_{\kk}(C^*[i], \kk) \; \cong \; \Hom_{\kk}( C^*, \kk) [i] $$
where following the usual convention for dual of (differential) graded objects  we define $$\Hom_{\kk}(C^*,\kk)_n := \Hom_{\kk}(C^n, \kk), \quad \Hom_{\kk}(D_*,\kk)^n := \Hom_{\kk}(D_n, \kk)$$ for any integer $n$.

\smallskip

 For a \emph{persistent} object $P$ (over $\Delta^+$ or $\R^2$, see~\ref{D:PerModule}) we will also use the standard notation $P[s]$ (where $s$ is in $\R^2$) for its shifted by the vector $s$, which is also a persistent object (Definition~\ref{D:shiftedpersitent}). Note that the shift is by a vector, i.e. a point in $\R^2$ not an integer.  We will also use letters such as $t$, $x$ or $\vec{\epsilon}$, $s_x$, $s_y$ for these operations.  This should cause no confusion since the sets of letters used in the two types of shifts are disjoint.
 
 \smallskip
 
 For instance for a graded persistent object $P^\bullet$, the notation $(P^\bullet[i])[s]=(P^\bullet[s])[i]$, where $i\in \Z$ and $s\in \R^2_{>0}$, stands for the  persistent object defined by $(P[i])[s] (x)^n:= P^{n+i}(x+s)$.

\section{The category of Mayer-Vietoris systems over $\R$}
In this section we study the notion of Mayer-Vietoris system which are persistence modules over $\{(x,y)\mid x+y > 0\}$ with additional structure. We first start by the latter notion. 
\subsection{Middle-exact persistence modules over $\Delta^+$}

We define $\Delta^+ := \{(x,y)\mid x+y > 0\}\subset \R^2$, equipped with the product (partial) order $(x,y)\leq (x',y')$ if and only if $x\leq x'$ and $y\leq y' $. 
In other words, $\Delta^+$ is the upper half-space above the antidiagonal $\Delta$ of $\R^2$, equipped with the induced product partial order of $(\R, \leq)\times (\R, \leq)$.

\smallskip

Recall that for $s=(s_1,s_2)\in \R^2_{>0}$, we denote $s_x:=(s_1,0)$ and $s_y := (0,s_2)$. 

\begin{defi}\label{D:PerModule}
A \emph{persistence module} over $\Delta^+$ is a functor $M : (\Delta^+,\leq) \longrightarrow \Mod(\kk)$, where $\Mod(\kk)$ is the category of $\kk$-vector spaces.

Persistence modules over $\Delta^+$ together with natural transformations of functors form a category, denoted $\Pers(\Delta^+)$.

Similarly the category of \emph{persistence comodules} over $\Delta^+$ is the category of functors $(\Delta^+,\leq)^{\op} \longrightarrow \Mod(\kk)$.
\end{defi}
In particular, the data of a  persistence module over $\Delta^+$ is encoded by the  structural  morphisms  
\begin{equation}\label{eq:structmorphism}
    M(v) \longrightarrow M(v+s)= M(v_1+s_1, v_2+s_2)
\end{equation}
defined for any $v=(v_1,v_2)\in \Delta^+$ and $s\in \R^2_{\geq0}\setminus\{0\}$ and their compatibilities. 

\smallskip

The motivating examples are the persistence modules $s=(s_1,s_2) \mapsto H_i(f^{-1}]-s_1,s_2[)$ associated, for any integer $i$, to the level-set of a continuous function $f: X\to \R$ (see example~\ref{Ex:MVfromFct}). Because of this traditional examples, we will sometimes use the terminology \emph{level-set persistence modules for arbitrary  persistence modules over $\Delta^+$ }.

\begin{defi}\label{D:shiftedpersitent} For $s\in \R^2_{\geq 0}\setminus\{0\}$ and $M\in \Obj( \Pers(\Delta^+))$, we define 
$M[s]$ to be the persistence module given, for any $v\in \Delta^+$, by $M[s](v)= M(v+s)$ with structural morphisms induced by those of $M$: for any $t\in \R^2_{\geq 0}\setminus \{0\}$:
$$ M[s](v)=M(v+s) \longrightarrow M(v+s+t)=M[s](v+t).$$
We extend the definition to $s=0$ by $M[0]=M$ with the identity for  $M\to M[0]$.
\end{defi}
It is immediate to check that the structural morphisms~\eqref{eq:structmorphism}
induces, for any $s\in \R^2_{\geq 0}\setminus\{0\}$ the \emph{canonical translation} maps:
\begin{equation}\label{eq:translmorphism}
   \tau_{s}^M: M \longrightarrow M[s].
\end{equation}
We will adopt the convention that when we do not need explicitly a notation for a translation or structural morphism we simply do not label it. 

\smallskip

For $M\in \Pers(\Delta^+)$, any $s\in\R^2_{>0} $ induces the short complex : 
\begin{equation}\label{eq:secforPersModule}
M\{s\} = M \longrightarrow M[s_x]\oplus M[s_y] \longrightarrow M[s] 
\end{equation}
where the first map is $\left( \begin{array}{c} \tau^M_{s_x}\\ -\tau^M_{s_y}\end{array}\right)$ and the second one is $(\tau^{M[s_x]}_{s_y}, \tau^{M[s_y]}_{s_x})$ in matrix notations. In other words for any $t\in \Delta^+$ and  $v\in M(t)$, the first map is given by $v\mapsto (\tau^M_{s_x}(v), - \tau^M_{s_y}(v))$ and the second by 
$M(t+s_x)\oplus M(t+s_y) \ni (v,w) \mapsto \tau^{M[s_x]}_{s_y}(v)+\tau^{M[s_y]}_{s_x}(w)$. 
The fact that~\eqref{eq:secforPersModule} is a complex is an immediate consequence of definition~\ref{D:PerModule}.

\begin{defi}\label{D:middleexact} An object
$M\in \Pers(\Delta^+)$ is said to be \textbf{middle-exact} if the complexes $M\{s\}$ are exact for every $s\in\R^2_{>0} $.
\end{defi}
\begin{remark}
We think of middle-exact complexes as being the analogue for the poset $\Delta^+$ of half the terms of the Mayer-Vietoris long exact sequence relating the various homology groups of two open subsets of a space, their reunion and intersection. What is missing to have a long exact sequence are precisely the connecting homomorphisms relating homology groups of different degrees. In Section~\ref{S:ClassMVSystems}, we will precisely introduce an additional data on a (graded) middle-exact object of $\Pers(\Delta^+)$ to obtain such long exact sequences. 
\end{remark}
Persistence modules have a barcode decomposition similar to peristence modules over $\R$ that we now describe. First we specify the various geometric types, called blocks, of the barcode. 
\begin{defi}\label{def:block_MV}
A block $B$ is a subset of $\R^2$ of the following type : 
\begin{enumerate}
    \item A \textbf{birthblock} (\textbf{bb} for short) if there exists $(a,b)\in \R^2$ such that $B = <a,\infty> \times <b,\infty>$, where $a$ and $b$ can eventually worth $-\infty$ simultaneously. Moreover, we will write that $B$ is of type \textbf{bb$^+$} if $a+b > 0$, and of type \textbf{bb$^-$} if $a+b \leq 0$. 
    \item A \textbf{deathblock} (\textbf{db} for short) if there exists $(a,b)\in \R^2$ such that $B = <-\infty,a> \times <-\infty,b>$. 
    Moreover, we will write that $B$ is of type \textbf{db}$^+$ if  $a+b>0$ and of type \textbf{db$^-$} if not. 
     \item A \textbf{horizontalblock} (\textbf{hb} for short) if there exists $a\in \R$ and $b\in \R\cup \{+\infty\}$ such that $B = \R \times <a,b>$.
     \item A \textbf{verticalblock} (\textbf{vb} for short) if there exists $a\in \R\cup \{+\infty\}$ and $b\in \R$ such that $B = <a,b> \times \R$.
\end{enumerate}
\end{defi}
\begin{remark}
Blocks are defined over the whole $\R^2$ and not just $\R^2_{>0}$.
\end{remark}

\begin{figure}
\begin{center}
\includegraphics[scale=0.9]{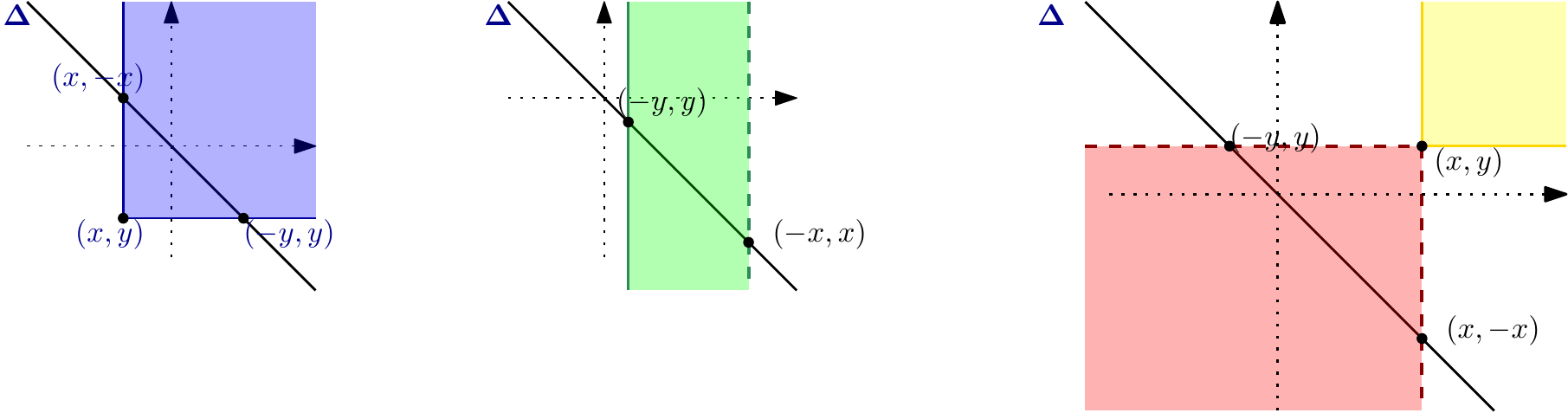}\label{fig:blocks}
\caption{On the left a block of type \textbf{bb}$^{-}$ pictured in blue. On the middle a block of type $\textbf{vb}$ pictured in green. On the right 
a block of type \textbf{db} in red and its dual block of type \textbf{bb}$^{+}$ in yellow.
The various coordinates refers to the intersection points of the boundaries of the blocks with the anti-diagonal $\Delta$ as well as the extremum of the birth or death blocks.
The dashed boundary lines means that the boundary line is not part of the block.}
\end{center}
\end{figure}
\begin{remark}\label{R;dualBlocks}
Note that a deathblock $B$ is characterized by its supremum\footnote{which is easily seen to be, if $B= <-\infty, x> \times <-\infty,y>$, the point $(x,y)\in \R^2$}, that is $\sup \{s\in B\}$ together with the data of whether its two boundary lines are in the block or not (note that the supremum is inside $B$ if and only if both boundaries lines are). 
Similarly a birth block $B'$  is characterized by its infimum $\inf \{s\in B'\}$ and the data of whether its boundary lines are in $B$ or not.  Note also that the vertical and horizontal blocks never have finite extremums.  
\end{remark}
\begin{defi}[Duality between death and birth blocks]
 \label{D:dualblocks} The dual  of a deathblock $B$ is the birthblock $B^\dagger$ whose infimum is the supremum of $B$ and whose vertical (resp. horizontal) boundary lines are in $B^\dagger$ if and only if the the vertical (resp. horizontal) boundary lines of $B$ are not. 
 
 Dually we define the dual $C^\dagger$ of a birthblock $C$ as the death block whose supremum is the infimum of $C$ and whose vertical (resp. horizontal) boundary lines are in $C^\dagger$ if and only if the the vertical (resp. horizontal) boundary lines of $C$ are not.
\end{defi}
\begin{remark}
 The rule $B\mapsto B^\dagger$ is involutive: $(B^\dagger)^\dagger=B$ and in particular exhibits a perfect duality between death and birth blocks.  Furthermore, note that the dual of a deathblock is of type \textbf{bb}$^{+}$ if and only if the deathblock has a non-trivial intersection with $\Delta^+$ i.e. is in \textbf{db}$^+$. 
\end{remark}

We now define the building blocks (that is the indecomposables in the middle exact case) of persistence modules over $\Delta^+$.
\begin{defi}\label{Def:blockmodule}
Let $B$ be a block, define the \textbf{block module} associated to $B$ by, for any $s\leq t \in \Delta^+$, : 
$$\kk^B(s) = \begin{cases} \kk \text{~if~}s\in B \\ 0 \text{~else}
\end{cases} ~~~\kk^B(s\leq t) = \begin{cases} \text{id}_\kk \text{~if~} (s,t)\in B^2 \\ 0 \text{~else} \end{cases} $$
\end{defi}
\begin{remark}\label{R:trivaildbmoins}
 If $B$ is in \textbf{db}$^{-}$, then $\kk^B=0$. Therefore we will usually not consider the block modules associated to such negative deathblocks. In what follows, the reader can safely assume that when we speak about a deathblock  we mean an element of \textbf{db}$^+$, unless otherwise stated.
\end{remark}

Let us denote, for $s\in \Delta^+$, $B -s=\{ t -s, \, t\in B\}$; this is a block of the same type as $B$.
\begin{lem}\label{L:shiftofBlock}
 Let $B$ be a block and $s\in \Delta^+$. There is a canonical isomorphism 
 $$\kk^B[s] \, \cong \, \kk^{B-s}. $$
\end{lem}
\begin{pf}
 By definition~\ref{Def:blockmodule}, we have that 
 $$\kk^B[s](t)= \kk^B(t+(s_1, s_2)) = \left\{ \begin{array}{ll} 
\kk & \text{if t }\in B -s \\
0 & \text{else.}\end{array}\right.$$
Therefore we have that $$\kk^B[s] \, \cong \, \kk^{B-s}.$$
\end{pf}
\begin{thm}[Cochoy-Oudot~\cite{CO17}, Botnan-Crawley-Boevey~\cite{BotCra18}] \label{thm:exactdecomp}
Let $M\in \Pers(\Delta^+)$ be middle exact and pointwise finite dimensional (pfd). Then there exists a unique multiset of blocks $\mathbb{B}(M)$ such that : 

$$M \simeq \bigoplus_{B\in \mathbb{B}(M)} \kk^B. $$
\end{thm}

\subsection{Mayer-Vietoris systems over $\R$ and their classification}\label{S:ClassMVSystems}

\begin{defi}\label{D:MVSystem}
We define the category $\text{M-V}(\R)$ of \textbf{Mayer-Vietoris persistent systems over} $\R$ as follows : 

\begin{itemize}
    \item[$\bullet$] Objects : collections $S=(S_i,\delta_i^s)_{i\in \Z,s\in \R_{>0}^2}$ where $S_i$ is in $\text{Pers}(\Delta^+)$ and $\delta_i^s\in \Hom_{\Delta^+}(S_i[s], S_{i-1})$, such that for all $i\in \Z$ and all $s\in \R_{>0}^2$, the following sequence 
    \begin{equation}\xymatrix{S_{i+1}[s] \ar[r]^-{\delta_{i+1}^s} & S_i \ar[r] & S_i[s_x] \oplus S_i[s_y] \ar[r] & S_i[s] \ar[r]^-{\delta_{i}^s} & S_{i-1}}\label{eq:MVSystemles}\end{equation}
    is exact and furthermore the following diagram  is commutative, for $s'\geq s$ :  \begin{equation}\label{D:ComSaquareforMV}\xymatrix{  S_i[s] \ar[r]^-{\delta_{i}^{s}} \ar[d] & S_{i-1} \ar[d]^{\text{id}_{S_{i-1}}} \\
    S_i[s'] \ar[r]^-{{\delta}_{i}^{s'}} & S_{i-1}.}\end{equation}
    
    \item[$\bullet$] Morphisms : for $(S_i,\delta_i^s)$ and $(T_i,\tilde{\delta}_i^s)$ two Mayer-Vietoris systems over $\R$, a morphism from $(S_i,\delta_i^s)$ to $(T_i,\tilde{\delta}_i^s)$ is a collection of morphisms $(\varphi_i)_{i\in \Z}$ where $\varphi_i \in \Hom_{\Pers(\Delta^+)}(S_i,T_i)$ such that the following diagram
    
   \begin{equation}\xymatrix{  S_i[s] \ar[r]^-{\delta_{i}^s} \ar[d]_{\varphi_i[s]} & S_{i-1} \ar[d]^{\varphi_{i-1}} \\
    T_i[s] \ar[r]^-{\tilde{\delta}_{i}^s} & T_{i-1}}\end{equation} commutes for all $i\in \Z$ and $s\in \R_{>0}^2$. 
\end{itemize}
For a Mayer-Vietoris system $S$ and $i\in \Z$, we will write $S_i$ for the associated object of $\Pers(\Delta^+)$ of $S$ which lies in degree $i$.
\end{defi}
A natural class of examples of such M-V systems is provided by homology of level-sets of a  continuous function on a topological space $X$. See, example~\ref{Ex:MVfromFct} below. Furthermore, we will see that any complex of sheaves $F^\bullet$ on $\R$ gives rise to a MV-system $\Psi(F^\bullet)$ (see Proposition~\ref{P:psiofCOnstrisMV}). 
\begin{remark}
\begin{itemize}
    \item Observe that if $(S_i,\delta_i^s)$ is a Mayer-Vietoris system, $S_i$ is in particular a middle exact modules, for $i\in \Z$. 
    \item The category $\text{M-V}(\R)$ is indeed a category. It is easy from the definition to observe that it is  additive. However, as we shall see later on, it is not abelian. 
\end{itemize}
\end{remark}

Our remaining goal in this section is to classify Mayer-Vietoris system in a way similar to Theorem~\ref{thm:exactdecomp}. For this, we introduce building blocks for those.
\begin{defi}\label{D:blocksmodulesforMV}
Let $\text{B}$ be a block (Definition~\ref{def:block_MV}) and $j\in \Z$. We define the Mayer-Vietoris system of degree $j$ associated to $\text{B}$, denoted $\text{S}^\text{B}_j$,   by : 

\begin{itemize}
    \item If B is of type \textbf{bb}$^-$, \textbf{hb} or \textbf{vb} then $\text{S}^\text{B}_j = (M_i,0)_{i,s}$ with $M_i = 0$ for all $i\not = j$ and $M_j = \kk_B$
    
    \item If B is of type \textbf{db}$^+$, then $\text{S}^\text{B}_j = (M_i,\delta_i^s)$ with $M_i = 0$ for all $i\not \in \{j+1,j\} $, $\delta_i^s=0$ for all $s\in \R^2$ and for $i\not = j+1$,  we define $M_{j+1} = \kk_{B^\dag}$, $M_{j}= \kk_{B}$, and $\delta_{j+1}^s : \kk_{B^\dag}[s] \to \kk_B$ by pointwise identities on $B^\dag [s]\cap B \cap \Delta^+$.

    \item Dually, if $B$ is of type \textbf{bb$^+$}, then define $S^{\text{B}}_j$ as $S^{B^\dag}_{j-1}$.
    \item If B is of type \textbf{db}$^{-}$, then we set $S^B=0$.
\end{itemize}  
\end{defi}
Of course the case of \textbf{db}$^{-}$ matches remark~\ref{R:trivaildbmoins}.

\begin{figure}
 \begin{center}
  \includegraphics{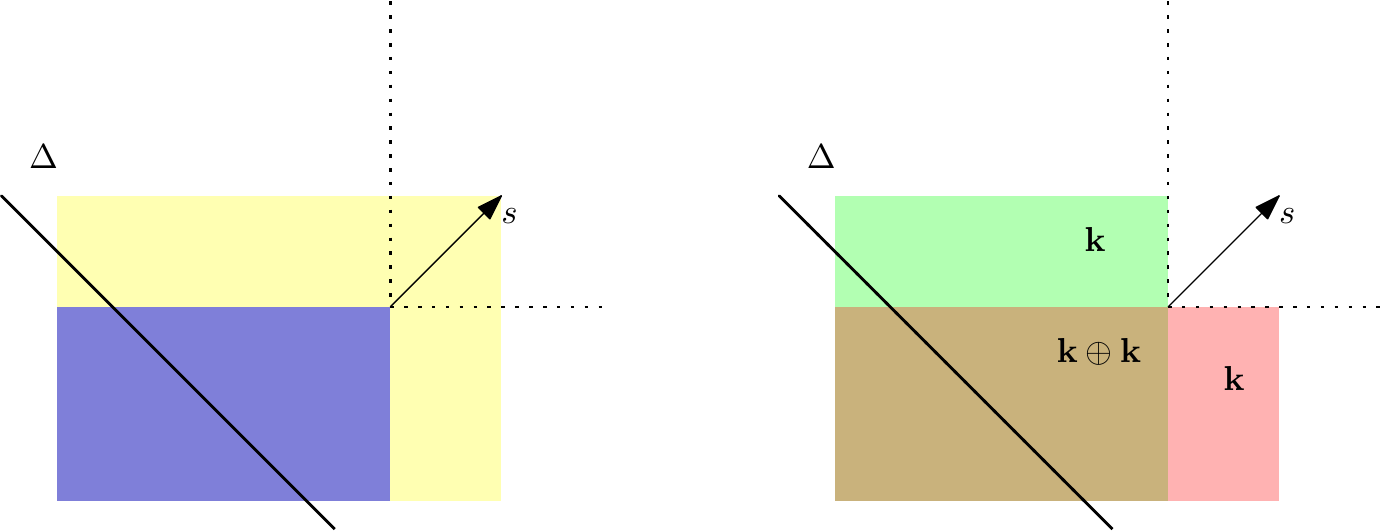}
  \caption{\emph{On the left:} A deathblock $B$ in yellow (and blue) and the translated deathblock $B-s=B[s]$ in blue where $s$ is the vector drawn. The dotted lines are the boundary of the dual birth block $B^{\dagger}[s]$. 
    \emph{On the right:} The  value of $\kk^B[s_x]\oplus \kk^B[s_y]$ in every region where the green region is the  translated death block $B-s_x$, the red is the translated death block $B-s_y$. The value is $0$ on the white region and below the antidiagonal.  }
  \label{fig:translationblocks}
 \end{center}
\end{figure}

\begin{remark}
One can easily see that the Mayer-Vietoris systems $S_j^B$ are indecomposable. We will refer to these Mayer-Vietoris systems to \emph{block MV-systems} for short.

Also note that for a block $B$ of type \textbf{db}$^+$ or \textbf{bb}$^{+}$ and $j\in \Z$, the graded persistent module $(M_i,0)_{i,s}$ with $M_i = 0$ for all $i\not = j$ and $M_j = \kk_B$ is \emph{not} a Mayer-Vietoris system.
\end{remark}
\begin{lem} The graded persistent modules $S_j^{B}$ associated to blocks $B$ in Definition~\ref{D:blocksmodulesforMV} are Mayer-Vietoris systems for any $j$ and block $B$.
\end{lem}
\begin{pf} We advise the reader to draw the different cases in a way similar to figure~\ref{fig:translationblocks}. Since $S_j^B =S_{j-1}^{B^\dagger}$, the case of \textbf{db}$^+$ and \textbf{bb}$^+$ are equivalent.

\smallskip 

 Note that every block which is not of type  \textbf{db} is stable by upward vertical and/or left-to-right  horizontal translations. It follows that $\kk^B \to \kk^B[s_x] \oplus \kk^B[s_y]$ is injective. Thus for blocks of type \textbf{vb}, \textbf{hb} or \textbf{bb}$^{-}$,  $S^B_j \to S^B_j[s_x]\oplus S^B_j [s_y]$ is one to one as well in every degree, a well as is the map $S^B_j \to S^B_j[s_x]\oplus S^B_j [s_y]$ in degree $i\neq j$ for $B$ of type \textbf{db} (and therefore also for $S^B_{j-1}$ if $B$ is of type \textbf{bb}$^{+}$ by definition~\ref{D:blocksmodulesforMV}). 

 \smallskip
 
 Note now that for a block $B$, if $z \in \R^2$ and $s\in \R^2_{>0}$ satisfies that $z+s \in B$, then either  $z+s_x$ or $z+s_y$ is in $B$ as well if $B$ is of type different from \textbf{bb}. Furthermore, for a block of type \textbf{bb}, the latter property only fails if $x \in B\cap B^\dagger$ where $B^\dagger$ is its dual (death)block. When $B$ is of type \textbf{bb}$^{-}$, those points are not in $\Delta^+$.   
 Therefore, the maps $\kk^B[s_x]\oplus \kk^B[s_y]\to \kk^B[s]$ are surjective for all blocks of type different from \textbf{bb}$^{+}$. 
 
 \smallskip 
 
 Let us now prove that the subsequences $ \xymatrix{\kk^B \ar[r] & \kk^B[s_x] \oplus \kk^B[s_y] \ar[r] & \kk^B[s]}$ are exact for any $B$; we have already seen that the composition is zero. Now, assume $(\alpha_x, \alpha_y) \in  \kk^B[s_x](v) \oplus \kk^B[s_y](v)$ is a \emph{nonzero} element in the kernel of $\tau_{s_x} \oplus \tau_{s_y}$. Then, if $v+s \in B$ then so are $v+s_x$ and $v+s_b$ and therefore $\tau_{s_x}$ and $\tau_{y}$ are the identity map $\kk \to \kk$. In particular $\alpha_x=\alpha_y=:\alpha$. But since $\kk^{B}[s](v)= \kk^{B}(s+v)=\kk$ as well, then $\kk^B(v)\to  
  \kk^B[s_x](v) \oplus \kk^B[s_y](v)$ is the map $(id, id)$ and hence $(\alpha, \alpha)$ is in its image. If $v+s\notin B$, then $B$ is not a birthblock and at least one element among $v+s_x$ and $v+s_y$ is not in $B$. If none are, then there is nothing to prove and if not then $B$ is either a vertical or horizontal block. In the first case, $v+s_y \in B$ and therefore $\kk^B \to \kk^B[s_y]$ is the identity map so that we have a preimage for $\alpha_y$. The other case is dual. This conclude the proof of the lemma for all blocks which are not of type \textbf{db}$^+$. 
 
 \smallskip 
 
 To prove the result for blocks of type \textbf{db}$^+$, since $S_i^B = S^B [-i]$ and by the injectivity result we have obtained at the beginning of that proof, it is enough to prove that the sequences 
 $$ \kk^{B^\dagger}[s](v) \to \kk^B(v) \to \kk^B[s_x](v)\oplus \kk^B[s_y](v)  $$
 are exact for any $s\in \R^2_{>0}$, $v\in \Delta^+$. If $v\notin B$, there is nothing to prove. Thus we assume $v\in B$.
 First, if both $\kk^B[s_x]$ and $\kk^B[s_y]$ are null, then  $x\in B\cap B^\dagger[s]= B\cap (B-s)^{\dagger}$. 
 Therefore, $\kk^{B^\dagger}[s](v) \to \kk^B(v)$ is the identity and the sequence is exact. 
 If both $\kk^B[s_x]$ and $\kk^B[s_y]$ are non-null, then $\kk^B(v) \to \kk^B[s_x](v)\oplus \kk^B[s_y](v) $ identifies with the necessarily injective diagonal inclusion and $v\notin B[s]^\dagger$ so that $\kk^{B^\dagger}[s](v)=0$ and the sequence is thus exact. 
 Finally if only one $\kk^B[s_x]$ or $\kk^B[s_y]$ is non-null, one of the map $\kk^B\to \kk^B[s_x]$ or $\kk^B\to \kk^B[s_y]$  is the identity-hence injective-and  we still have $v\notin B[s]^\dagger$. Thus  $\kk^{B^\dagger}[s](v)=0$. The sequence is again exact and the lemma is proved.
 \end{pf}

Denote by \textbf{M-V$^+ (\R)$} the full \emph{sub-category of} Mayer-Vietoris systems over $\R$ whose objects are the MV systems $S =(S_j,\delta_j^s)$ such that there exists $N\in\Z$ with $S_j = 0$ for all $j< N$. In other words, M-V$^+ (\R)$ is the subcategory of lower-bounded Mayer-Vietoris systems.

\begin{thm}[Classification of pfd M-V systems]\label{thm:decom_MV}

Let $S$ an object of M-V$^+ (\R)$ which is pointwise finite dimensional.  Then there exists a unique collection of multisets of blocks $\mathbb{B}(S) = (\mathbb{B}_j(S))_{j\in \Z}$ of type \textbf{bb$^-$}, \textbf{hb}, \textbf{vb}, and \textbf{db}$^+$, such that we have an isomorphism in M-V$(\R)$ : 

$$S \simeq \bigoplus_{j\in \Z} \bigoplus_{B\in \mathbb{B}_j(S)} S_j^B $$

We call $\mathbb{B}(S)$ the \textbf{barcode} of $S$. It completely determines $S$ up to isomorphism of Mayer-Vietoris systems.
\end{thm}

\begin{remark}
By Definition~\ref{def:block_MV}, birth blocks of type \textbf{bb$^+$} generate the same MV systems as their dual death blocks, therefore  they come in pairs in the decomposition given by Theorem~\ref{thm:decom_MV}, which explains why the blocks of type \textbf{bb$^+$} are ignored in the barcode.
\end{remark}

To prove the theorem~\ref{thm:decom_MV}, we will use the following technical lemmas :

\begin{lem}\label{lem:dbnull}
Let $S = (S_j,\delta_j^s)$ be a pfd MV-system over $\R$. If $\mathbb{B}(S_j)$ contains only blocks of type \textbf{db}$^+$,
then $S=0$.

\end{lem}

\begin{pf}
Given $s\in \R^2_{>0}$, 
the universal property of cokernels and the exactness of~\eqref{eq:secforPersModule} imply that $\delta^s_j$ factorizes through $$\text{coker}\left (S_j[s_x] \oplus S_j[s_y] \longrightarrow S_j[s] \right ).$$
Now, this cokernel is trivial 
since by assumption $S_j$ is isomorphic to a direct sum of blocks of type \textbf{db}$^+$. Therefore, $\delta^s_j = 0$.

Consequently, for every $B\in \mathbb{B}(S_j)$ and every $s\in \R^2_{s>0}$, the exact sequence of persistence modules 
$$0 \longrightarrow \kk_B \longrightarrow \kk_B[s_x] \oplus \kk_B[s_y] $$
yields $B = \emptyset$ since $B$ is assumed to be of type \textbf{db}$^+$.
\end{pf}

\begin{lem}\label{lem:splitbb+}

Let $S$ be a pfd MV-system over $\R$, such that there exists  a block $B = \langle a,\infty \rangle \times \langle b, \infty \rangle$ of type \textbf{bb$^+$} such that $B\in\mathbb{B}(S_{j})$ for some $j\in\Z$. Then there exists a pfd MV system $\Sigma$ such that : 

$$ S \simeq S_{j}^{B} \oplus \Sigma = S_{j-1}^{B^\dag} \oplus \Sigma  $$
\end{lem}

\begin{pf}
Let $s \geq \sqrt{2}(a+b,a+b)$, since $\text{coker}\left (\kk_B[s_x] \oplus \kk_B[s_y] \longrightarrow \kk_B[s] \right ) \simeq \kk_{B^\dag} $  we have the following commutative diagram, where the rows are exact sequences and where $\varphi$ exists (and is injective) by the universal property of cokernels: 
 
 \begin{equation}\label{eq:MV-summand}
 \begin{gathered}
 \xymatrix{0 \ar[r] & \kk_B \ar@{^{(}->}[d] \ar[r] &\kk_B[s_x] \oplus \kk_B[s_y] \ar[r] \ar@{^{(}->}[d] & \kk_B[s]\ar[r] \ar@{^{(}->}[d] & \kk_{B^\dag} \ar[r] \ar@{.>}[d]^\varphi  & 0 \\
\dots \ar[r] & S_j \ar[r] & S_j[s_x] \oplus S_j[s_y] \ar[r] & S_j[s] \ar[r]^{\delta_j^s} & S_{j-1} \ar[r]& \dots  }
\end{gathered}
\end{equation}
 
Since $B^\dag$ is a directed ideal of $\Delta^+$, $\kk_{B^\dag}$ is an injective object of $\Pe(\Delta^+)$ by lemma 2.1 of\cite{BotCra18}. Therefore, $\varphi$ splits and $\text{im} \varphi \simeq \kk_{B^\dag}$ is a summand of $S_{j+1}$.  The commutativity of~(\ref{eq:MV-summand}) then implies the existence of a complement $X_j$ of $\text{im} (\kk_B\hookrightarrow S_j)$ in $S_j$, such that $S$ decomposes locally as follows:

$$
\xymatrix{X_j \oplus \kk_B \ar^-{\simeq}[d]\ar[r] &X_j[s_x] \oplus X_j[s_y] \oplus  \kk_b[s_x] \oplus \kk_b[s_y] \ar^-{\simeq}[d]\ar[r] & X_j[s] \oplus \kk_B[s] \ar^-{\simeq}[d] \ar[r]  & X_{j-1} \oplus \kk_{B^\dag} \ar^-{\simeq}[d]\\
S_j \ar^-{\sigma}[r] & S_j[s_x] \oplus S_j[s_y] \ar[r] & S_j[s] \ar[r]^{\delta_j^s} & S_{j-1}
}$$

Note that we may assume without loss of generality that $X_j \supseteq \ker \sigma$.
Then, by exactness of~$S$, we have $\text{im} \delta_{j+1}^s = \ker \sigma\subseteq X_j$, therefore our local decomposition extends to a full decomposition of $S$, which means that the upper row complex in~(\ref{eq:MV-summand}) is a summand of~$S$.
\end{pf}

\begin{pf}[of theorem~\ref{thm:decom_MV}] 
Let $S =(S_j,\delta_j^s) \in $ M-V$^+ (\R)$, and assume without loss of generality that the lower bound~$N$ is equal to~$1$. Then, all the $S_j$'s are middle-exact pfd persistence modules over $\Delta^+$, therefore they decompose uniquely (up to isomorphism) as direct sums of block modules, by theorem \ref{thm:exactdecomp}. Note that for $j\leq 0$ the decomposition is trivial.

\medskip

\noindent {\bf The finite barcode case:}

We first show the result in the case where $\mathbb{B}(S_j)$ is finite for every $j\in\Z$.  For each $j\in\Z$, fix an isomorphism $\varphi_j : S_j \stackrel{\sim}{\longrightarrow}\bigoplus_{B\in\mathbb{B}(S_j)}\kk_B$. Thus, the family $(\varphi_j)_j$ induces an isomorphism of MV systems from $S$ to 
    $$S' := \left(\bigoplus_{B\in\mathbb{B}(S_j)}\kk_B, ~\varphi_{j-1}\circ \delta_j^s \circ \varphi_j^{-1}[s] \right )_{j\in\Z,s\in \R_{>0}^2}$$ 
  
  Let $B\in \mathbb{B}(S_j)$ of type either \textbf{bb}$^-$, \textbf{hb} or \textbf{vb}, then for $s\in \R_{>0}^2$, the map :  $$\kk_B[s_x]\oplus\kk_B[s_y]\longrightarrow \kk_B[s]$$ is surjective. Thus,  $\varphi_{j-1}\circ \delta_j^s \circ \varphi_j^{-1}[s]$ is zero on $\kk_B[s]$. This proves that $S^B_j$ is a summand of $S'$. Finally, noting $\mathbb{B}^{-}(S_j)$ the multi-set of intervals of $\mathbb{B}(S_j)$ of type either \textbf{bb}$^-$, \textbf{hb} or \textbf{vb}, we have : 
  
 $$S' = \left (\bigoplus_{j\in\Z} \bigoplus_{B\in  \mathbb{B}^{-}(S_j) } S^B_j \right ) \oplus \left (\bigoplus_{j\in\Z} \bigoplus_{B\in\mathbb{B}(S_j)\backslash \mathbb{B}^{-}(S_j)}\kk_B, ~\varphi_{j-1}\circ \delta_j^s \circ \varphi_j^{-1}[s] \right ).$$ 
 
 There remains to prove that the right-hand side of the direct sum, noted $S''$, decomposes in MV($\R$). For $j\in\Z$, the barcode $\mathbb{B}(S_j)\backslash \mathbb{B}^{-}(S_j)$ contains only blocks of type either $\textbf{bb$^+$}$ or $\textbf{db}^+$. Denote by $\mathbb{B}(S_j)^+$ the multiset of blocks of type \textbf{bb$^+$} involved in $\mathbb{B}(S_j)$.
 Let us prove by induction that, for any  $j_0\geq 1 (=N)$, there exists a MV system $\Sigma^{j_0}$ such that
 
 $$S'' \simeq \left ( \bigoplus_{1\leq j < j_0} \bigoplus_{B\in \mathbb{B}(S_j)^+ } S_j^B \right ) \oplus \Sigma^{j_0}.$$
 
 For $j_0 = 1$ the property clearly holds with $\Sigma^{j_0} = S''$.
 Let us now assume the property holds up to some $j_0\geq 1$. Since  $\mathbb{B}(S_{j_0})^+$ has finite cardinality, Lemma~\ref{lem:splitbb+} (applied repeatedly) decomposes $\Sigma^{j_0}$ as
 
 $$\Sigma^{j_0} \simeq \left(\bigoplus_{B\in \mathbb{B}(S_{j_0})^+} S_{j_0}^B\right) \oplus \Sigma^{j_0+1},$$
 which yields the induction step.

 Now, given $j_0\geq 1$, for any $j< j_0$ the barcode of $\Sigma_{j}^{j_0}$ can only contain deathblocks by construction. Therefore, by Lemma~\ref{lem:dbnull}, we have $\Sigma_{j}^{j_0} = 0$.
It follows that
 $$S" \simeq \bigoplus_{j\in\Z} \bigoplus_{B\in \mathbb{B}(S_j)^+} S^B_j,$$
thus concluding the decomposition in the finite barcode case.
 
\medskip

\noindent {\bf The infinite barcode case:}

We now generalize to the case where the barcodes $\mathbb{B}(S_j)$ can be infinite. For the same reason as in the finite case, each block of type \textbf{bb}$^-$, \textbf{hb} or \textbf{vb} involved in some barcode $\mathbb{B}(S_j)$ splits as a summand $S_j^B$ of $S$. Hence, we are reduced to proving the existence of the decomposition in the case where $S$ is a pfd MV system, and $\mathbb{B}(S_j)$ contains only blocks of type \textbf{bb$^+$} or \textbf{db}$^+$ for all $j\in\Z$. Given $n\in \Z_{>0}$, define $\Delta^+_n := \Delta^+\cap\{(x,y)\in \R^2 \mid  x\leq n , y\leq n\}$. Define also 
$$\mathbb{B}(S_j)_n := \{B\in \mathbb{B}(S_j) \mid B ~\text{is of type \textbf{bb$^+$} and } B\cap \Delta^+_n \not = \emptyset ~ \text{or~} B \text{~is of type \textbf{db}$^+$ and } B \subset \Delta^+_n \}.$$

Then it is clear that $\mathbb{B}(S_j) = \bigcup_n \mathbb{B}(S_j)_n$, and since $S$ is pointwise finite dimensional, $\mathbb{B}(S_j)_n$ contains finitely many blocks of type \textbf{bb$^+$}, for all $n \geq 0$.
We now identify each $S_j$ with its block decomposition via some fixed isomorphism, and for $n\geq 0$ we define ${}_n\tilde{S}$ as follows: 

$${}_n\tilde{S} = \left (\bigoplus_{B\in \mathbb{B}(S_j)_n } \kk_{B}, (\delta_j^s)_{|\bigoplus_{B\in \mathbb{B}(S_j)_n } \kk_{B}} \right )  $$

 Let us prove that $\tilde{S}$ is a sub-MV system of $S$. To do so, it is sufficient to prove that for all $j\in  \Z$, the image of $(\delta_j^s)_{|\bigoplus_{B\in \mathbb{B}(S_j)_n}}$ is contained in $\bigoplus_{B\in \mathbb{B}(S_{j-1})_n } \kk_{B}$. Fix $j\in \Z$ and $s\in \R^2_{>0}$. Then $(\delta_j^s)_{|\bigoplus_{B\in \mathbb{B}(S_j)_n}}$ factorizes uniquely through: 
 
 \begin{align*}
     & \text{coker} \left ( \bigoplus_{B\in \mathbb{B}(S_j)_n } \kk_{B}[s_x] \oplus  \kk_{B} [s_y] \longrightarrow \bigoplus_{B\in \mathbb{B}(S_j)_n } \kk_{B}[s]  \right )   \\ 
     &\simeq \bigoplus_{B\in \mathbb{B}(S_j)_n }  \text{coker} \left ( \kk_{B}[s_x] \oplus  \kk_{B} [s_y] \longrightarrow  \kk_{B}[s] \right ) \\
     &=  \bigoplus_{ \substack{B\in \mathbb{B}(S_j)_n\\ B ~\text{is of type \textbf{bb$^+$}} } } \text{coker} \left ( \kk_{B}[s_x] \oplus  \kk_{B} [s_y] \longrightarrow  \kk_{B}[s] \right ) 
 \end{align*} 
 
As previously, for every $B\in \mathbb{B}(S_j)_n$ of type \textbf{bb$^+$}, we can find $s\in \R^2_{>0}$ such that the canonical map: 

$$ \text{coker} \big ( \kk_{B}[s_x] \oplus  \kk_{B} [s_y] \longrightarrow  \kk_{B}[s] \big  )  \longrightarrow \bigoplus_{B\in \mathbb{B}(S_{j-1}) } \kk_{B} $$

is a monomorphism. And as seen in the proof of Lemma~\ref{lem:splitbb+}, $ \text{coker} \big ( \kk_{B}[s_x] \oplus  \kk_{B} [s_y] \longrightarrow  \kk_{B}[s] \big  ) $ is isomorphic to $\kk_{B^\dag}$, hence an injective object of $\Pers(\Delta^+)$, so its image splits off as a summand of $\bigoplus_{B\in \mathbb{B}(S_{j-1}) } \kk_{B}$ and is therefore included in  
$$  \kk_{B^\dag}^m\subset \bigoplus_{B\in \mathbb{B}(S_{j-1}) } \kk_B $$

where $m$ is the multiplicity of $B^\dag$ in $\mathbb{B}(S_{j-1}) $. Since $B^\dag\in \mathbb{B}(S_{j-1})_n$, we conclude that $\text{im} ((\delta_j^s)_{|\bigoplus_{B\in \mathbb{B}(S_j)_n}}) \subset \bigoplus_{B\in \mathbb{B}(S_{j-1})_n } \kk_{B} $. This proves that ${}_n\tilde{S}$ is a sub-MV system of $S$.

Then, we can apply our decomposition result in the finite barcode case to ${}_n\tilde{S}$. And since we have the filtration

$$S = \bigcup_{n\geq 0} {}_n\tilde{S} $$

which stabilizes pointwise, we get a decomposition for $S$.
\end{pf}
\begin{remark}
 Let us finish by a remark on the \lq\lq{}derived\rq\rq{} meaning of  Mayer-Vietoris systems. The axioms and structure we put on  $\text{M-V}(\R)$ are actually encoding a natural homotopy property. To state it, we have to consider the (derived) category of \emph{2-parameter persistence chains complexes}, that is the (associated derived) category of functors $\Delta^+ \to \text{dg-}\Mod(\kk)$.   Taking the direct sum of homology groups of 2-parameter persistence chain complex gives a graded 2-parameter persistence module. 
 Such graded 2-parameter persistence modules $(H_i(C_\bullet)_{i\in \Z}$ that can be lifted to a Mayer-Vietoris system are precisely those such that the underlying 2-parameter persistence chain complex $C_\bullet$ satisfies the following property\footnote{which is best expressed using $\infty$-categories or model categories}: 
 \begin{itemize}
     \item[ ] For any $s\in \Delta^+$, the canonical map $ C_\bullet[s_x] \oplus C_\bullet[s_y] \to C_\bullet [s]$ exhibits $ C_\bullet [s]$ as the \emph{homotopy} quotient
 $\mathbf{hocoker} \Big(C_\bullet \to C_\bullet[s_x] \oplus C_\bullet[s_y]\Big)$ of  the persistence chain complex morphisms
 $C_\bullet \to C_\bullet[s_x] \oplus C_\bullet[s_y]$.
 \end{itemize}
 A down to earth way of expressing this homotopy quotient property is to say $C_\bullet [s]$ is quasi-isomorphic  to the cone of 
 $C_\bullet \to C_\bullet[s_x] \oplus C_\bullet[s_y]$ as a persistent chain complex over $\Delta^+$. In other words, the structure of Mayer-Vietoris systems is essentially encoding the data of a homotopy property carried by their underlying chain complexes; property expressing that the chain complex at a $(x, y)+s$ is  determined by those of the chain complexes at the point $(x,y)$, $(x,y) +s_x$, $(x,y)+s_y$ for any $s\in \Delta^+$ which exhibits a local to global coherence of the values of those special 2-parameter persistence chain complexes.
\end{remark}

\subsection{Interleaving distance for M-V systems}
We have a (fully faithful) functor $(\R ,\leq) \to (\R^2, \leq)$ given by the diagonal embedding  $t\mapsto \vec{t}$ where $\vec{t}=(t,t)$. \emph{We} also \emph{denote} $\overrightarrow{(-)}: (\R_{>0} ,\leq) \to (\Delta^+, \leq)$  the induced functor. 

\smallskip

Given $\varepsilon \geq 0$, and $M = (M_i,\delta_i^s)$ a Mayer-Vietoris system over $\R$ (as in Definition~\ref{D:MVSystem}), observe that the collection $\tau_{\vec{\varepsilon}}^M:=(\tau_{\vec{\varepsilon}}^{M_i})_{i\in\Z}$ is a morphism of Mayer-Vietoris systems$M \longrightarrow M[\vec{\varepsilon}]$, where $M[\vec{\varepsilon}]:=(M_i[\vec{\varepsilon}],\delta_i^s[\vec{\varepsilon}])$.

\begin{defi}
Let $M$ and $N$ two Mayer-Vietoris systems over $\R$. An \textbf{$\varepsilon$-interleaving} between $M$ and $N$ is the data of two morphisms of M-V systems $f = (f_i) : (M_i,\delta_i^s) \longrightarrow (N_i[\vec{\varepsilon}],\tilde{\delta}_i^s [\vec{\varepsilon}])$ and $g = (g_i) : (N_i,\tilde{\delta}_i^s) \longrightarrow (M_i[\vec{\varepsilon}],\delta_i^s[\vec{\varepsilon}]) $ such that the following diagram commutes : 

\begin{equation}\label{eq:Interleaving} \xymatrix{
M  \ar[rd]\ar@/^0.7cm/[rr]^{\tau_{2\vec{\varepsilon}}^M} \ar[r]^{f} & N[\vec{\varepsilon}]  \ar[rd] \ar[r]^{g[\vec{\varepsilon}]} & M[2\vec{\varepsilon}] \\
N  \ar[ur]\ar@/_0.7cm/[rr]_{\tau_{2\vec{\varepsilon}}^N} \ar[r]^{g} & M[\vec{\varepsilon}]  \ar[ur] \ar[r]^{f[\vec{\varepsilon}]} & N[2\vec{\varepsilon}]
   } \end{equation}
   
If $M$ and $N$ are $\varepsilon$-interleaved, we shall write $M\sim_\varepsilon^{MV} N$

\end{defi}

\begin{defi}\label{D:InterleavingforMV}
Define the interleaving distance between two Mayer-Vietoris systems $M$ and $N$ to be the non-negative or possibly infinite number : 

$$d_I^{MV}(M,N) := \inf \{\varepsilon \geq 0 \mid M \sim_\varepsilon^{MV} N  \}. $$
\end{defi}
\begin{remark}\label{R:DMVgeneralizedDstandard}
 The interleaving distance for Mayer-Vietoris persistence system is just a derived\footnote{because we precisely requires the morphisms to commute with the maps $\delta_s^i$ connecting homology groups of different degrees. This claim will be even more supported by the isometry theorem~\ref{T:MainIsometry}} extension of the usual interleaving distance in $\Pe(\Delta^+)$ defined, for $M, N \in \Obj(\Pe(\Delta^+))$ by 
 $$d_I(M,N) := \inf \{\varepsilon \geq 0 \mid M \sim_\varepsilon^{\Delta^+} N  \} $$
 where $M \sim_\varepsilon^{\Delta^+} N$ means that $M$ and $N$ are $\varepsilon$-interleaved as persistence modules, that is there  exists $f: M\to N[\vec{\varepsilon}]$ and $g: N\to M[\vec{\varepsilon}]$ are persistence modules morphisms satisfying that the diagram~\eqref{eq:Interleaving} commutes. 
\end{remark}

 We say that a Mayer-Vietoris system $M=(M_i,\delta_i)_{i\in \Z}$ is \textbf{bounded} if there is only finitely many $M_i$ which are non-zero.

We now turn to a main source of examples of Mayer-Vietoris systems. 
\begin{ex}[Mayer-Vietoris system associated to continuous functions]\label{Ex:MVfromFct}
Let $h: X\to \R$ be a continuous function on a topological space $X$. For any $
x=(x_1,x_2)\in \Delta^+$, we 
 set $$M_i^h(x):= H_i(h^{-1}(]-x_1, x_2[).$$
 If $x'=(x'_1, x'_2)\geq x$, then we have the inclusion $]-x_1, x_2[\subset ]-x_1', x_2'[$ inducing, for all $i$'s, homomorphisms $M_i^h(x)=H_i(h^{-1}(]-x_1, x_2[) \to H_i(h^{-1}(]-x'_1, x'_2[)= M^h(x')$ in homology. 
 By Lemma~\ref{L:Defiota}, this makes $M_i^h(-)$ a persistence module over $\Delta^+$, which is called the \textbf{level set persistence module associated to $h:X\to \R$}. 
 
 \smallskip
 
 Now, let $s=(s_1,s_2) \in \R^2_{>0}$. For any $x=(x_1,x_2)\in \Delta^+$, we have that the open interval $]-x_1-s_1, x_2,+s_2[ $ has a cover given by the two open sub-intervals $]-x_1-s_1, x_2[$ and $]-x_1, x_2,+s_2[ $ whose intersection is $]-x_1, x_2[$. Therefore the Mayer-Vietoris sequence associated to this cover gives us linear maps $(\delta_{i}^{s,x})_{i\in \N}$ and exact sequences 
  \begin{equation}\xymatrix{M^h_{i+1}[s](x) \ar[r]^-{\delta_{i+1}^{s,x}} & M^h_i(x) \ar[r] & M^h_i[s_x](x) \oplus M^h_i[s_y](x) \ar[r] & M^h_i[s](x) \ar[r]^-{\delta_{i}^{s,x}} & M^h_{i-1}}(x).\label{eq:MVSystemfromFct}\end{equation}
  We write $\delta^s_i: M_i^h[s] \to M_{i-1}^h$ the maps given at every point $x$ by $\delta_i^{s,x}$ and for $i\leq 0$ we set $\delta^i=0$.
  \begin{prop}\label{P:MVfromFct}
  The  $\delta^s_i$'s are persistence modules morphisms and  makes the collection $M^h:=(M^h_i, \delta^s_i)_{i\in \Z, s\in \R^2_{>0}}$ a Mayer-Vietoris persistence system over $\R$.
  
  Furthermore, the assignment $f\mapsto M^f$ is a functor $M^{(-)}: \mathrm{Top}_{|\R} \longrightarrow \text{M-V}(\R)$.
  \end{prop}
  \begin{pf}
  The fact that the $\delta^s_i$ are persistence modules maps as well as the commutativity of diagram~\eqref{D:ComSaquareforMV} follow from the naturality of the Mayer-Vietoris sequence. The exactness of~\eqref{eq:MVSystemfromFct} implies the condition~\eqref{eq:MVSystemles}. 
  
  Recall from the introduction that $\mathrm{Top}_{|\R}$ is the category of topological spaces over $\R$ which by definition has objects given by continuous functions $f:X\to \R$ where $X$ is a topological space. The set of morphisms from    $f:X\to \R$ to $g: Y\to \R$ is the set of all continuous maps $\phi:X\to Y$ such that the diagram $\xymatrix{X\ar[r]_{\phi} \ar@/^1pc/[rr]^{f} & Y\ar[r]_{g} & \R } $ is commutative. Since $f^{-1}(]-x,y[) =\phi^{-1} \big( g^{-1}(]-x,y[)\big)$, we have that $\phi$ restricts to a continuous map $\phi: f^{-1}(]-x,y[) \hookrightarrow g^{-1}(]-x,y[)$. Therefore we have  induced maps $\phi_*(x,y):M_i^f=H_i(f^{-1}(]-x,y[) \to H_i(g^{-1}(]-x,y[)=M_i^g$ after taking homology for all $(x,y)\in \Delta^+$. The functoriality of the homology functor and Mayer-Vietoris sequence prove that this $\phi_*$ is a morphism of Mayer-Vietoris system and furthermore that the assignment $f\mapsto M^f$, $\phi\mapsto (s\mapsto\phi_*(s))$ is a functor. 
  \end{pf}
\end{ex}
\begin{ex}
 Assume $X$ is a smooth or topological manifold and $h: X\to \R$ is continuous. Then the Mayer Vietoris system $M_i^h$ (given by example~\ref{Ex:MVfromFct}) is bounded since an open subset of a manifold is a manifold and hence has no homology in degrees higher than its dimension.  
\end{ex}

In particular, we obtain from the degree-wise stability of interleaving distance between level-set persistence modules the following : 

\begin{prop}\label{P:StabilityfordMV}
Let $h_1,h_2 : X \to \R$ two continuous functions defined on the topological space $X$. Then : $$d_I^{MV}(M^{h_1},M^{h_2})\leq \sup_{x\in X} |h_1(x)-h_2(x)| $$

\end{prop}
\begin{pf} If the distance is $\infty$, there is nothing to prove. Otherwise, let $\varepsilon= \sup_{x\in X} |h_1(x)-h_2(x)|$. Then for any $(x,y)\in \Delta^{+}$, we have  level-set inclusions  $h^{-1}_1(]-x,y[) \subset h^{-1}_2(]-x-\varepsilon, y+\varepsilon[)$ and  $h^{-1}_2(]-x,y[) \subset h^{-1}_1(]-x-\varepsilon, y+\varepsilon[)$ which induce persistence modules over $\Delta^+$ morphisms 
$$f:\big(M_i^{h_1}(x,y)=H_i(h_1^{-1}(]-x,y[) \to H_i(h_2^{-1}(]-x-\varepsilon,y+\varepsilon[)= M_i^{h_2}[\vec{\varepsilon}](x,y)\big)_{(x,y)\in \Delta^+},  $$ 
$$g:\big(M_i^{h_2}(x,y)=H_i(h_2^{-1}(]-x,y[) \to H_i(h_1^{-1}(]-x-\varepsilon,y+\varepsilon[)= M_i^{h_1}[\vec{\varepsilon}](x,y)\big)_{(x,y)\in \Delta^+} $$ since 
taking homology groups is a functor and by lemma~\ref{L:Defiota}. 

The fact that these maps are Mayer-Vietoris systems morphisms follows again as in proposition~\ref{P:MVfromFct} by the naturality of the Mayer-Vietoris sequence associated to open covers of the intervals $]-x-\varepsilon, y+\varepsilon[$ by $]-x-\varepsilon, y[$ and $]-x, y+\varepsilon[$.
\end{pf}

\section{Stable sheaf theoretic interpretation of persistence}
The relationship between (co)sheaf theory and persistence homology has been emphasized by the work of Curry~\cite{Curr14} and Kashiwara-Schapira~\cite{KS18, Kash18}. We will recall basics of this point of view in this section and construct a functor from (level-set) persistent objects (over $\Delta^+$) to sheaves.

\smallskip

In this paper, we follow the standard notations of \cite{Kash90}, \cite{KS18}, \cite{Berk18} for sheaves. 

In particular,  $\kk$ will denote a field, $\Mod(\kk)$ the category of vector spaces over $\kk$ and, for   a topological space $X$,  we will note $\Mod(\kk_X)$ the category of sheaves of $\kk$-vector spaces on $X$ and $\text{PSh(X)}$ the category of presheaves of $\kk$-modules on $X$. For shortness, we will also write $\Hom$ for $\Hom_{\Mod(\kk_\R)}$.

\smallskip

Henceforth, $\D^b(\kk)$ will be the \emph{derived} category of complexes of $\kk$-modules with bounded cohomology, and $\D^b(\kk_X)$ will be the one of complexes of sheaves with bounded cohomology of $\kk$-modules over $X$. Recall that the derived category is obtained from $\Mod(\kk_X)$ by inverting quasi-isomorphisms of complexes of abelian sheaves.  
 Unless the context is unclear, we will simply use the word  sheaf for an object of $\D^b(\kk_X)$. 
We will use the standard Grothendieck operations on sheaves as in~\cite{Kash90}.  
Note that, associated to any open-closed subset $Z$ of $X$, is a sheaf $\kk_Z$ whose main property is that its stalks are $\kk$ at any point in $Z$ and are $0$ else, see~\cite{Kash90}.
 
\subsection{Convolution distance for sheaves after Kashiwara-Schapira}
In~\cite{KS18} Kashiwara and Schapira have defined a (pseudo)distance on the derived category of sheaves. This distance is a sheaf version (derived by design)  of the interleaving distance of persistence modules. It is based on convolution of sheaves which we now explain. 

Let  $\V$  be an euclidean vector space, which in our case of interest will simply be $\R$. We let $s: \V\times \V \to \V$ be the addition map $(t,t')\mapsto t+t'$.  

\begin{defi}[\cite{KS18}]\label{D:Convolution} The convolution of sheaves  $\D^b(\kk_\V)\times \D^b(\kk_\V) \to \D^b(\kk_\V)$ is the bifunctor given, 
for $F,G\in \Obj(\D^b(\kk_\V))$,  by the formula : $$F\star G = R s_!(F\boxtimes G)$$ where $\boxtimes$ is the external tensor product of sheaves. 
\end{defi}
To define the convolution distance,  we will only need a very specific case :  convolution by the constant sheaf supported on a ball centered at 0. 
More precisely we define, for $\varepsilon \in \R$, 
\begin{equation}\label{eq:defKepsilon}
 K_\varepsilon :=  \left\{ \begin{array}{cc} 
 \kk_{\{x\in \V,  \|x\|\leq \varepsilon\}}  & \mbox{if } \varepsilon \geq 0 \\
 \kk_{\{x\in \V, \|x\|< -\varepsilon\}}[\dim(\V)]  & \mbox{if } \varepsilon < 0 .\end{array} \right.
\end{equation}
The convolution by $K_\varepsilon$ has some nice properties : 
\begin{prop}[\cite{KS18}]\label{P:propertiesofconvolution} Let $\varepsilon, \varepsilon'\in \R$ and $F \in \Obj(\D^b(\kk_\V))$. 
\begin{enumerate}
\item There are functorial isomorphisms $F\star K_0\simeq F$ and $(F\star K_{\varepsilon} )\star K_{\varepsilon'} \simeq F \star K_{\varepsilon + \varepsilon'} $.
\item If $\varepsilon' \geq \varepsilon $, there is a canonical morphism of sheaves 
$K_{\varepsilon'}\to K_{\varepsilon}$ in $\D^b(\kk_\V)$
inducing a natural transformation $F\star K_{\varepsilon'} \to F \star K_{\varepsilon} $. 
In the special case where $\varepsilon = 0$, we simply write $\phi_{F, \varepsilon'}$ for this natural transformation.
\end{enumerate}
\end{prop}
Convolution by $K_{\varepsilon}$ is the sheaf analogue of the canonical shift $F[\varepsilon]$ of a persistence module. Indeed, by
Proposition~\ref{P:propertiesofconvolution}.(1),  for any map $f: F\star K_{\varepsilon} \to G$ we get  canonical maps  
\begin{equation} \label{eq:propertiesofconvolution} f\star K_{\tau}: F\star K_{\varepsilon+\tau}\simeq F\star  K_{\varepsilon}\star K_{\tau}  \to G \star K_{\tau}.\end{equation} 
These  maps allow us to define interleaving.
\begin{defi}\label{D:inteleavingconvolution}
\begin{enumerate}\item For $F,G\in \Obj(\D^b(\kk_\V))$ and $\varepsilon\geq 0$, one says that $F$ and $G$
are $\varepsilon$-\textbf{interleaved} if there exists two morphisms in $\D^b(\kk_\V)$, $f : F\star K_\varepsilon \to G$ 
and $g :  G \star K_\varepsilon \to F$ such that the compositions
$F\star K_{2\varepsilon} \stackrel{ f\star K_\varepsilon}{\longrightarrow}K_\varepsilon\star G \stackrel{g}{\longrightarrow} F $ and $G\star K_{2\varepsilon} \stackrel{g \star K_\varepsilon}{\longrightarrow}K_\varepsilon\star F \stackrel{f}{\longrightarrow} G $ are the natural morphisms $F\star K_{2\varepsilon} \stackrel{\phi_{F,2\varepsilon}}{\longrightarrow} F$ and $G\star K_{2\varepsilon} \stackrel{\phi_{G,2\varepsilon}}{\longrightarrow} G$, that is, we have a commutative diagram in $\D^b(\kk_\V)$ :

$$ \xymatrix{
F\star K_{2\varepsilon}  \ar[rd]\ar@/^0.7cm/[rr]^{\phi_{F,2\varepsilon}} \ar[r]^{f\star K_\varepsilon} & G\star K_\varepsilon  \ar[rd] \ar[r]^{g} & F\\
G \star K_{2\varepsilon}  \ar[ur]\ar@/_0.7cm/[rr]_{\phi_{G,2\varepsilon}} \ar[r]^{g\star K_\varepsilon} & F\star K_\varepsilon  \ar[ur] \ar[r]^{f} & G
   } $$ 
In this case, we write $F \sim_\varepsilon G$.
\item or $F,G\in \Obj(\D^b(\kk_\V))$, we define their \textbf{convolution distance} as : $$d_C(F,G) : = \text{inf} \left (\{ +\infty \} \cup \{a \in \R_{\geq 0 } \mid \text{$F$ and $G$ are $a$-isomorphic} \} \right )$$ \end{enumerate}
\end{defi}
The convolution distance has the following properties
\begin{prop}[\cite{KS18} and \cite{Berk18}]\label{P:StabforConv}\begin{enumerate}
\item The convolution distance is a closed extended pseudo-metric on $\D^b(\kk_\V)$ that is, for $F,G,H \in \Obj(\D^b(\kk_\V))$ : 
\begin{enumerate}
\item $d_C(F,G) = d_C(G,F)$,
\item $d_C(F,G) \leq d_C(F,H) + d_C(H,G)$,
\item  if moreover $F$ and $G$ are constructible (see definition below), one has $d_C(F,G) \leq \varepsilon 
\iff F \sim_\varepsilon G $.
\end{enumerate}
\item \textbf{(Stability Theorem)} If $X$ is a locally compact topological space, and $f,g : X \to \V$ are continuous functions, then for any $F\in \Obj(\D^b(\kk_\V))$ one has : $$d_C(\Rr f_* F, \Rr g_* F) \leq \sup_{x\in X} \|f(x) - g(x) \|$$
and the same is true for $\Rr f_!$ and $\Rr g_!$.
\end{enumerate}
\end{prop}

\subsection{Graded barcodes and derived isometry theorem}
There is a notion of barcodes for \emph{constructible} sheaves that mimicks the persistence case. This allows  to define a derived bottleneck distance following~\cite{Berk18}.

\begin{defi}
A sheaf  $F\in \Obj(\Mod(\kk_\V))$, $F$ is said to be \textbf{constructible} if there exists a locally finite stratification of $\V = \sqcup_\alpha S_\alpha$, such that for each stratum $S_\alpha$ is locally closed in $\V$, the restriction $F_{|S_\alpha}$ is locally constant and furthermore,  the stalks $F_x$ are of finite dimension for every $x\in \V$. 

We  write respectively $\Mod_{\R c}(\kk_M)$ and $\D^b_{\R c}(\kk_M)\cong \D^b(\Mod_{\R c}(\kk_M))$ for the  category of constructible sheaves on $M$ and  the full (triangulated) subcategory of $\D^b(\kk_M)$ consisting of complexes of sheaves whose cohomology objects lies in $\Mod_{\R c}(\kk_M)$.
\end{defi}
Note that the notion of constructibility is precisely what is usually called $\R$-constructibility. Since no other notion will show up in this work we simply drop the $\R$. 
\begin{remark} The condition on the stalks is the sheaf analogue of the condition of being pointwise finite dimensional  for persistence modules.   
\end{remark}

There is a decomposition similar to persistence for constructible sheaves. Namely we have the following two results.
\begin{thm}[Decomposition - \cite{KS18} Theorem 1.17.]\label{T:KSdecomposition}
Let $F \in \Obj(\Mod_{\R c}(\kk_\R))$, then there exists a locally finite family of intervals $\{I_\alpha\}_{\alpha \in A}$ such that $F \simeq \mathop{\bigoplus}\limits_{\alpha \in A} \kk_{I_\alpha}$. Moreover, this decomposition is unique up to isomorphism.
\end{thm}
\begin{cor}[Structure]\label{T:KSstructure}
Let $G^\bullet\in \Obj(\D^b_{\R c}(\kk_\R))$. 
\begin{enumerate} \item Then there exists an isomorphism in $\D^b_{\R c}(\kk_\R)$ : $$G^\bullet \simeq \bigoplus_{j\in\Z} \Ho^j(G^\bullet)[-j]$$
where $\Ho^j(G^\bullet)$ is seen as a complex of sheaves concentrated in degree $0$.
\item For each $j\in \Z$, there is a unique multiset $\B^j(G^\bullet) $ of intervals such that  $\Ho^j(G^\bullet)\simeq \bigoplus_{I \in \B^j(G^\bullet)} \kk_{I}$. 
\end{enumerate}\end{cor}
This corollary allows us to define the graded barcode of an object of $\D^b_{\R c}(\kk_\R)$ following~\cite{Berk18} which is the  derived enhancement of the usual barcode of 1d-persistence modules.
\begin{defi}\label{D:gradedbarcode}
 The \textbf{graded-barcode} $\B^\bullet(G^\bullet)$ of $G^\bullet$ is the sequence of multisets $(\B^j(G^\bullet))_{j\in \Z}$.
 
We write $\B_C^\bullet(G^\bullet)$, $\B_L^\bullet(G^\bullet)$ and $\B_R^\bullet(G^\bullet)$ for the  sub-multisets of $\B^\bullet(G^\bullet)$ consisting respectively of the closed or bounded open intervals, semi-open intervals which are open on the right,  semi-open intervals which are open on the right (and not equal to $\R$). 

The interval appearing in the respective subsets will be called respectiveley of central type, left type and right type.
 \end{defi}
By the corollary~\ref{T:KSstructure}, the graded-bracode uniquely determines the complex of sheaves up to isomorphisms in $\D^b_{\R c}(\kk_\R)$ and furthermore we have a unique decomposition $$G^\bullet \cong G_C^\bullet \oplus G_L^\bullet \oplus G_R^\bullet $$ into sheaves whose cohomology only have supports    
in intervals of central type, left type and right type respectively. This is called the CLR decomposition in~\cite{Berk18}.
\begin{lem}[\cite{Berk18}]\label{L:Caracinterleaving} Let $F^\bullet, \, G^\bullet \in \D^b_{\R c}(\kk_\R)$.
If $\varepsilon \in \R$, then $G^\bullet$ and $F^\bullet$ are $\varepsilon$-interleaved if and only if 
 $G_C^\bullet \mathop{\sim}\limits_{\varepsilon} F_C^\bullet$, $G_L^\bullet \mathop{\sim}\limits_{\varepsilon} F_L^\bullet$ and $G_R^\bullet \mathop{\sim}\limits_{\varepsilon} F_R^\bullet$. In particular, 
  bars of a given type can only be interleaved with bars of the same type.
\end{lem}
One can characterize the geometric  condition for two sheaves $\kk_I$ and $\kk_J$ on intervals of  same types to be $\varepsilon$-interleaved in terms of their endpoints. 
We will need the following result.
\begin{prop}[Proposition 3.8 in \cite{Berk18}]\label{P:P38Berk18}
Let $\varepsilon \geq 0$, and $a\leq b$ in $\R \cup \{\pm \infty\}$. Then : 
\begin{itemize}
\item[$\bullet$] $\kk_{[a,b]}\star K_\varepsilon \simeq \kk_{[a-\varepsilon, b + \varepsilon]} $,
\item[$\bullet$] $\kk_{]a,b[}\star K_\varepsilon \simeq 
\begin{cases}
\kk_{]a+\varepsilon,b-\varepsilon[} ~~\text{if} ~~\varepsilon < \frac{b-a }{2} \\
\kk_{[b- \varepsilon , a+ \varepsilon]}[-1] ~~\text{if} ~~\varepsilon \geq \frac{b-a }{2}\, ,
\end{cases}$

\item[$\bullet$] $\kk_{]a,b]}\star K_\varepsilon \simeq \kk_{]a+\varepsilon, b + \varepsilon]} $,
\item[$\bullet$] $\kk_{[a,b[}\star K_\varepsilon \simeq \kk_{[a-\varepsilon, b - \varepsilon[} $.
\end{itemize}
\end{prop}

 In order to define the bottleneck distance, we first define the notion of $\varepsilon$-matching. 
 
\begin{defi}\label{D:Epsilonmatching}
Let $\B$ and $\B'$ be two graded-barcodes and $\varepsilon \geq 0$. An $\varepsilon$-\textbf{matching} between $\B$ and $\B$ is the data of
\begin{enumerate}
\item partial matchings: $\sigma_R^j: \B^j_R \not \to \B_R^j$, $\sigma_L^j: \B^j_L \not \to (\B')_L^j$ for all $j\in \Z$
satisfying that 
\begin{description}
\item $(i)$ for any matched pair $I$, $\sigma_R(I)$  (resp. $J$, $\sigma_L (J)$), 
one has $\kk_I[-i] \sim_\varepsilon\kk_{\sigma_R^j(I)}[-i]$ (resp. $\kk_J[-j] \sim_\varepsilon\kk_{\sigma_L^j(J)}[-j]$) and
\item $(ii)$ for  the $I \in \B_R$ and $\B_L$ which are not matched,  one has $\kk_I[-i]\sim_\varepsilon 0$.\end{description}
\item a \emph{bijection} $\sigma^j_C :  \B_C \longrightarrow \B'_C$ satisfying,  for any $I \in B^j_C $, that  $\kk_I \sim_\varepsilon \kk_{\sigma^j_C(I)}[-\delta]$ and
\begin{description}
\item $(i)$ that $\delta =0 $ if $I$ and $\sigma^j_C(I)$ are both open or both closed, 
\item $(ii)$ and $\delta = 1$ if $I$ is open and $\sigma^j_C(I)$ is closed, $\delta = -1$ if  $I$ is closed and $\sigma^j_C(I)$ is open.
 \end{description}
\end{enumerate}
\end{defi}
\begin{remark}
 An $\varepsilon$-matching can match bars of different degrees, but only if one of them is compact and  the other one is open and differs in degree by $+1$. In all other situations one can only match bars of same degrees. See~\cite{Berk18}.
\end{remark}

\begin{defi}
Let $\B$ and $\B'$ be two graded-barcodes, then one defines their \textbf{bottleneck distance} to be the possibly infinite positive value: $$d_B(\B,\B') = \text{inf} \{\varepsilon \geq 0 \mid \text{there exists a } \varepsilon  \text{-matching between } \B \text{ and } \B' \} $$
\end{defi}
The graded bottleneck distance is isometric to the convolution.
\begin{thm}[Isometry~\cite{Berk18}]\label{T:DerivedIsometry}
Let $F^\bullet,\, G^\bullet$ be two objects of $\D^b_{\R c}(\kk_\R)$. Then 
$$ d_C(F,G) =  d_B(\B(F),\B(G))$$
\end{thm}

\subsection{Extending level-set persistence modules as pre-sheaves over $\R$}\label{SS:levelsettopresheaves}
In this section we interpret level-set persistence modules as (pre)sheaves on the line $\R$. 
 Let $(\Ouv(\R), \subset )$ be the poset of open subsets of $\R$ ordered by the inclusion. We denote in the same way the associated category.  
\begin{lem} \label{L:Defiota}Set $\iota: (\Delta^+, \leq)\to (\Ouv(\R), \subset)$ to be given on objects by 
  $\iota: s=(s_1, s_2) \mapsto ]-s_1, s_2[$.
  Then $\iota$ is a well defined  fully faithfull functor.
  The esential image of $\iota$ is precisely the full subcategory of bounded open intervals of $\R$. 
\end{lem}
In particular, restricting to objects of those categories,  $\iota$ is a bijection from $\Delta^+$ to bounded open intervals of $\R$. 
\begin{pf}
 By definition  $$ s=(s_1, s_2)\in \Delta^+ \;\Longleftrightarrow \;-s_1<s_2$$ hence $\iota$ is well defined, injective on objects with image the bounded open intervals.  Furthermore, if $(s_1, s_2)\leq (s_1', s_2')$ then 
 $-s_1'\leq -s_1<s_2\leq s_2'$ which proves that $\iota$ is order presearving (and necessarily fully faithfull since the morphisms are empty or a singleton). 
\end{pf}

Given $M \in \Obj(\Pe(\Delta^+))$ we can consider its pointwise dual $ t\mapsto \Hom_{\Mod(\kk)}(M(t);\kk)$ which has a canonical structure of a persistence comodule, that is of an object of  
$$\text{Fun}((\Delta^{+})^{\op}; \Mod(\kk))\,\cong\, \text{Fun}(\Delta^+; \Mod(\kk)^{\op})^{\op}. $$  We denote by $M^*$ this dual of $M$. More precisely, $M^*$ is the composition of functors
$$M^* :=\, \Delta^{+\op} \stackrel{M^{\op}}\longrightarrow \Mod(\kk)^{\op} \stackrel{\Hom_{\Mod(\kk)}(-;\kk)}\longrightarrow \Mod(\kk). $$

Since $M^*$ is a persistence comodule, for any  open $U \subset \R$, we have a module 
\begin{equation}\label{eq:defMtilde} 
\Tilde{M} (U):=  \varprojlim_{]-x,y[\subset U} M^*((x,y)).\end{equation}
                                     
\begin{lem}\label{L:valueoftilde} There is a functor $\tilde{(-)}: \Pe(\Delta^+)\to \text{PSh}(\R)^{\op}$   extending the formula~\eqref{eq:defMtilde} into a canonical presheaf on $\R$, that is such that for   $U \in \Obj(\Ouv(\R))$, one has $$\Tilde{M} (U):=  \varprojlim_{]-x,y[\subset U} M^*((x,y)).$$ 
\end{lem}
\begin{pf} One notice that the formula exhibits  $\Tilde{M}$ as a   Kan extension which makes it into a presheaf canonically. Indeed, 
 consider $\iota^\text{op} : (\Delta^{+})^{\op}\longrightarrow \Ouv(\R)^{\op}$ the (opposite of the) functor defined previously (see~\ref{L:Defiota}) and let  $\text{Ran}_{\iota^\text{op}} M^*$  be the right Kan extension along $\iota^\text{op}$ of $M^*$, which is therefore by definition an object of $\text{PSh}(X)$ : 
$$\xymatrix{ \Delta^{+\text{op}}
\ar@{^{(}->}[r]^-{\iota^{\text{op}}} \ar[d]_{M^*} & \Ouv(\R)^\text{op} \ar@{.>}[dl]^{\text{Ran}_{\iota^\text{op}} M^* =: \Tilde{M}} \\
\Mod(\kk) & } $$
As $\Mod(\kk)$ is complete,  the pointwise formula~\eqref{eq:defMtilde} is an immediate consequence.  
\end{pf}
\begin{remark}
 A corollary of the proof is that the restriction morphism of $\Tilde{M}$ are given, for $U\subset V$, by the canonical restrictions $M^*((a,b)) \to M^*((x,y))$ for any $ U\supset ]-x,y[ \subset ]-a,b[\subset V$ and the induced (by the universal property) map on the limits. 
\end{remark}
Composing $\tilde{(-)}$ with the (opposite of the) sheafification functor $ \text{PSh}(\R) \to \Mod(\kk_{\R})$ gives the functor from persistence modules on $\Delta^{+}$ to sheaves on $\R$. 
\begin{defi}\label{D:defofBar}
We set $\overline{M}$ to be the sheafification of the presheaf $\tilde{M}$ and we write 
$\Xi:=\overline{(-)} : \Pe(\Delta^+)\longrightarrow \Mod (\kk_\R)^{\op}$ for the induced functor $M\mapsto \overline{M}$. We call $\overline{(-)}$ the level-set persistence to sheaves functor. 
\end{defi}
Similarly there is a functor going in the other direction defined as follows. Given a sheaf on $\R$, by restriction to open intervals and using the identification of lemma~\ref{L:Defiota}, we get a persistence comodule. Since pointwise duality transforms a persistence comodule into a persistence module we obtain the functor 
\begin{equation}
 \label{eq:defpi}  \pi : \Mod(\kk_\R)^{\op} \longrightarrow \Pe(\Delta^+), \quad F\mapsto \Hom_{\Mod(\kk)}(F_{ | \text{open intervals}};\kk)
\end{equation}
where $F_{ | \text{open intervals}}$ is the composition $F\circ \iota^{op}$ (we recall that the image of $\iota^{op}$ is precisely the open subintervals).
\emph{We call} (abusively) this functor the \emph{restriction-to-intervals functor}. 

We will also write $\text{bidual}_{\Mod(\kk_\R)}$ the endofunctor of  sheaves  which to a sheaf $M$ 
associate its pointwise bidual $N\mapsto (N^*)^*$. There is a canonical natural transformation 
\begin{equation}\label{eq:unitadjPetoSh} \text{id}_{\Mod(\kk_\R)} \longrightarrow \text{bidual}_{\Mod(\kk_\R)}\end{equation} given by the pointwise canonical morphism.
\begin{prop}\label{P:PropertiesofBar}
The level-set persistence to sheaves functor $\Xi : \Pe(\Delta^+)\longrightarrow \Mod (\kk_\R)^{\op}$ from Definition~\ref{D:defofBar} satisfies the following properties.
\begin{enumerate}
    \item Its composition  with the restriction-to-intervals functors is the canonical biduality functor: $\Xi\circ \pi = \text{bidual}_{\Mod(\kk_\R)}$. 
    
    In particular the restriction of this composition of functors to the subcategory of pointwise finite dimensional\footnote{where we mean the sheaves whose stalk at each point are finite dimensional} objects is naturally isomorphic to $\text{id}_{\Mod(\kk_\R)^{pfd}}$
    \item It is the right adjoint of the functor $\pi$: 
    $$\Hom_{\Pe(\Delta^+)}(M, \pi(F)) \; \cong \; \Hom_{\Mod(\kk_\R)^{\op}}(\overline{M}, F)\; = \; \Hom_{\Mod(\kk_\R)}(F,\overline{M}) $$ and the map~\eqref{eq:unitadjPetoSh} is the counit of this adjunction.
    \item If $M$ is a pointwise finite dimensional persistent module, and $M \simeq \oplus_i M_i$, then $\overline{M} \simeq \oplus_i \overline{M_i} $
    \item Assume that $M \in \Obj(\Pe(\Delta^{+}))$ is pointwise finite dimensional. Then, for all $\alpha \in \R$, we have natural isomorphisms  $$ \varprojlim_{]-x;y[\ni \alpha}M((x,y)) \;\simeq\; \tilde{M}_{\alpha}\; \simeq \; \overline{M}_{\alpha}$$ provided that the left hand side is finite dimensional.
    
    \item One can identify $\overline{M}$ with the image of  pre-sheaves morphism : $\Tilde{M}\longrightarrow \prod_{\alpha\in\R}\Tilde{M}_{\alpha} $.
\end{enumerate}
\end{prop}

\begin{pf}
\begin{enumerate}
    \item First, once the formula  $\Xi\circ \pi= \overline{(-)}\circ \pi =\text{bidual}_{\Mod(\kk_\R)}$ is proved, to check  that the asserted restriction of the composite $\overline{(-)}\circ \pi $ is canonically isomorphic to the identity, it is sufficient  to prove that the canonical transformation $\text{id}_{\Mod(\kk_\R)} \to \text{bidual}_{\Mod(\kk_\R)}$ is an isomorphism on all stalks when restricted to a pointwise finite dimensional sheaf. This reduces the statement to the standard case of finite dimensional vector spaces.  To prove the formula note that the value of a sheaf  $F$ is a  on an open $U$ is uniquely determined by its value on any open cover; furthermore, for any open interval $I$ one has that $$\varprojlim_{]-x,y[\subset I} ((F(I))^*)^* \cong (F(I)^*)^*.$$ In particular, on can restrict to cover by open intervals $(I_j)$ of $U$ and compute the value of the sheaf at $U$ as a limit.
   Therefore, noticing that the intersection of two intervvals is an interval, we get that for a sheaf $F$, one has 
   \begin{multline*}\overline{\pi(F)}(U) = \varprojlim \Big(\xymatrix{\overline{\pi(F)}(\coprod U_k) \ar@<1ex>[r] \ar@<-1ex>[r] & \overline{\pi(F)}(\coprod (U_i \cap U_j)) } \Big) \\ 
   \cong \; \varprojlim \Big(\xymatrix{\prod (F(U_k)^*)^* \ar@<1ex>[r] \ar@<-1ex>[r] & \prod(F(U_i \cap U_j)^*)^* } \Big) \\ \cong \; (F(U)^*)^*
   \end{multline*}
   \item Since the sheafification functor is a right adjoint and by universal property of Kan extensions  we have natural isomorphisms:
   \begin{multline*}
  \Hom_{\Mod(\kk_\R)^{\op}}(\overline{M}, F)  = \Hom_{\Mod(\kk_\R)}(F,\overline{M})    \\ \cong \;
  \Hom_{\text{PSh}(\R)}(F, \tilde{M}) \; \cong \; \Hom_{\text{PSh}(\R)}(F, \text{Ran}_{\iota^{\op}} M^*) \\ \cong \; \Hom_{\Pe\big((\Delta^{+})^{ op}\big)}(F_{ | \text{op. intervals} }, \Hom_{\kk}(M, \kk)) \\ \cong \; \Hom_{\kk}(F_{ | \text{op. intervals} }\mathop{\otimes}\limits_{\Delta^+} M, \kk) \\ \cong \; 
  \Hom_{\Pe(\Delta^+)}(M, \Hom_{\kk}(F_{ | \text{op. intervals} }, \kk)) \\ \cong \; 
 \Hom_{\Pe(\Delta^+)}(M, \pi(F)). 
   \end{multline*}
   Here for a persistent comodule $L$ and a persistent module $M$, we denote $L\mathop{\otimes}\limits_{\Delta^+} M$ the coequalizer $$\mathrm{coeq}\left(\coprod_{s_1\stackrel{\phi}\to s_2\in \Delta^+}  L(s_2) \otimes_{\kk} M(s_1) \stackrel{\phi^*\otimes id}{\underset{id\otimes \phi_*}{\rightrightarrows}} \coprod_{s\in \Delta^+} L(s)\otimes_{\kk} M(s) \right)$$ where the upper and lower maps are induced by the (co)persistence structures. 
   
This proves the adjunction formula. The fact that the counit is given by~\eqref{eq:unitadjPetoSh} is 
a direct consequence of  the proof of property (1).
    \item The functors $\Hom(-;\kk)$, sheafification (which is a right adjoint) as well as right Kan extensions  commute with finite direct sums. This gives the finite sums case. But the assumption ensures it is enough to estblish the resutls on the stalks and therefore the canonical map $\bigoplus \overline{M}_i \to \overline{\bigoplus M_i}$ is an isomorphism.
    \item Write $\text{Int}(\alpha)$ for the (full) subcategory of $\Ouv(\R)$ consisting of intervals containing $\alpha$. Let us fix  $G :  (]0,\infty[,\leq)  \longrightarrow \text{Int}(\alpha)$ defined by 
  $G(\varepsilon) = ]\alpha - \varepsilon , \alpha + \varepsilon[  $. Then $G$ is a functor and is initial among functors $ (]0,\infty[,\leq)  \longrightarrow \text{Int}(\alpha)$. Therefore : 
     $$ \varprojlim_{]-x;y[\ni \alpha}M((x,y)) \simeq \varprojlim M\circ G =\varprojlim_{\varepsilon > 0}M((\varepsilon - \alpha  , \alpha + \varepsilon)) $$

     Since $(]0,\infty[,\leq) $ is a totally ordered set, we can apply the theorem of decomposition of pfd modules over totally ordered sets to $M\circ G$, thus there exists a multiset $\mathbb{B}(M\circ G)$ of intervals of $\R$ such that : 
     \begin{equation}\label{eq:DefMoG} M\circ G \simeq \bigoplus_{I\in \mathbb{B}(M\circ G)} \kk_I \end{equation}
    It follows that \begin{equation}\label{eq:porjlimM} \varprojlim_{\varepsilon > 0}M((\varepsilon - \alpha  , \alpha + \varepsilon)) \simeq \prod_{\substack{I\in \mathbb{B}(M\circ G)\\ 0 \in \text{closure}(I)}} \kk  \end{equation}

     Now if $\varprojlim\limits_{]-x;y[\ni \alpha}M((x,y))$ is finite dimensional then the above product in the right hand side of~\eqref{eq:porjlimM} is a finite product and thus a direct sum: 
     $\prod_{\substack{I\in \mathbb{B}(M\circ G)\\ 0 \in \text{closure}(I)}} \kk \; \simeq \; \bigoplus_{\substack{I\in \mathbb{B}(M\circ G)\\ 0 \in \text{closure}(I)}} \kk.$ Therefore we have 
     \begin{align*}
        \varprojlim_{]-x;y[\ni \alpha}M((x,y))
        & \simeq \bigoplus_{\substack{I\in \mathbb{B}(M\circ G)\\ 0 \in \text{closure}(I)}} \kk  \\
        &\simeq \varprojlim_{]-x;y[\ni \alpha} \Hom \left ( M((x,y)), \kk \right ) \; \text{ (by~\eqref{eq:DefMoG} and finite dimensionality)}\\ 
        &\simeq \tilde{M}_\alpha \; \text{ (by~\eqref{eq:defMtilde})}\\
        &\simeq \overline{M}_\alpha.
    \end{align*}
    \item This is a general fact for sheaves on a $T_1$-topological space, that is for sheaves on a space for which all   points are closed.
\end{enumerate}
\end{pf}
\begin{remark} We have sticked in this paper to the traditionnal point of view of looking at level-set as being given by homology functors and thus as persistent objects; point of view for which computational models are well developped. This is the reason why some (bi)duality shows up in the picture. It is possible (and actually slightly easier) to construct an analogue of $\Xi: M\mapsto \overline{(M)}$ going from persistence comodules to sheaves.
\end{remark}

Let $\Delta = \{(-x,x)\mid x\in\R\}$, and $p : \Delta \longrightarrow \R$ be the projection $(x_1,x_2)\mapsto x_2$ onto the second coordinate. Recall that for any block $B$ (Definition~\ref{def:block_MV}) we have defined (see~\ref{Def:blockmodule}) a persistence module $\kk^B\in \Obj(\Pe(\Delta^+))$.
\begin{prop}\label{P:BarofBlock}
Let $B$ be a block. Let $a,b\in \R$ be such that $<a,b> = p(B\cap \Delta)$, with the convention that $a=1$ and $b=-1$ when $ p(B\cap \Delta) = \emptyset$.  
\begin{enumerate}
    \item If $B$ is of type \textbf{dquad}, then $\overline{\kk^B} \simeq \kk_{]a,b[}$.
    \item If $B$ is of type \textbf{bquad}, then $\overline{\kk^B} \simeq \kk_{[a,b]}$.
    \item If $B$ is of type \textbf{vquad}, then $\overline{\kk^B} \simeq \kk_{]a,b]}$.
    \item If $B$ is of type \textbf{hquad}, then $\overline{\kk^B} \simeq \kk_{[a,b[}$.
\end{enumerate}
\end{prop}
\begin{pf}
Let $B$ be of type \textbf{dquad}. If $B$ is included in $\mathbb{R}^2\setminus \Delta^+$, then $\kk^B$ is identically null as well as $\kk_{]a,b[}$ and there is nothing to prove. If not, $B$ has a non-trivial intersection with $\Delta$ and   $(-a,b)=\sup_{B}(s)$ are the coordinates of the supremum of $B$ for the order relation of $\Delta^+$. Then, for  $s=(s_1,s_2)\in \mathbb{R}^2_{>0}$, one has 
$$\kk^B(s) =\left\{ \begin{array}{ll} \kk & \mbox{if } (s_1,s_2)< (-a,b) \\
0 & \mbox{ if } (s_1, s_2) \not\leqslant (-a,b)\end{array}\right . $$ 
Hence $\kk^B(s)$ is non-zero  if $ a\leq -s_1 <s_2 \leq b$ and always null if either $-s_1<a$ or $s_2 >b$. It follows that 
that for $\alpha \in \mathbb{R}$, then $ \kk^B(-\alpha-\varepsilon, \alpha+\varepsilon) =0$ 
if $\alpha \notin ]a,b[$ and for all $\alpha \in ]a,b[$, there exists $\eta >0$ such that 
$\kk(-\alpha -\eta, \alpha+\eta) =\kk$. We conclude that
\begin{equation}\label{eq:MVofdbquad}\varprojlim_{]-x, y[ \ni \alpha} \kk^B(x,y) = \left\{ \begin{array}{ll} \kk & \mbox{if } \alpha \in ]a,b[ \\
0 & \mbox{ else. }\end{array}\right . 
 \end{equation}
By claim 5 of Proposition~\ref{P:PropertiesofBar}, we deduce that $\overline{\kk^B} \cong \kk_{]a,b[}$.
Similarly, 
 if $B_v$ is a vertical block, delimited by the lines $x=-b$, $x=-a$ with $a<b$, then we have
$$\kk^{B_v}(s_1, s_2) =\left\{ \begin{array}{ll} \kk & \mbox{if } -b<s_1< -a \\
0 & \mbox{ if } s_1 >-a \mbox{ or } s_1<-b \end{array}\right . $$
independently of whether the boundary lines are part of $B_v$ or not. In particular, for any  $\alpha \in ]a,b]$, there exists $\eta >0$ such that 
$\kk(-\alpha -\eta, \alpha+\eta) =\kk$ while there exists $\varepsilon >0$ such that
$ \kk^B(-\alpha-\varepsilon, \alpha+\varepsilon) =0$ if $\alpha \leq a$ or $\alpha >b$. 
As in the \textbf{dquad} case~\eqref{eq:MVofdbquad}, we thus find that 
\begin{equation*}\varprojlim_{]-x, y[ \ni \alpha} \kk^{B_v}(x,y) = \left\{ \begin{array}{ll} \kk & \mbox{if } \alpha \in ]a,b] \\
0 & \mbox{ else. }\end{array}\right . 
 \end{equation*}
 The last two other cases are obtained using a similar analysis. 
\end{pf}

\begin{remark}
In particular, $\overline{\kk^B}$ does not depend on whether $B$ contains its boundary or not. If $B$ is of type $\textbf{bb}^+$, then $\overline{\kk^B} =0$ (since $B\cap \Delta =\emptyset$). 
\end{remark}
\begin{remark}[characterizations of $(a,b)$]\label{R:caracofab}
 If $B$ is of type $\textbf{bb}^-$, then the numbers $a$ and $b$ are characterized by the fact that the point $(-b,a)$ is the infimum of the points in $B$, see figure~\eqref{fig:blocks}.
 
 Similarly, if $B$ is of type $\textbf{dquad}$, then the numbers $a$ and $b$ satisfies that the point $(-a,b)$ is the supremum of the points in $B$. 
 
 Finally, for $B$ of type \textbf{vquad}, $a$ and $b$ satisfies that $B$ has boundary given by the lines of equation $x=-b$ and $x=-a$, while if it is of type \textbf{hquad}, $a$ and $b$ satisfies that the boundary of $B$ are the horizontal lines of equations $y=a$ and $y=b$.
\end{remark}

Blocks of type \textbf{db}$^{+}$, \textbf{hb}, \textbf{vb} and  \textbf{bb}$^{-}$ are actually uniquely determined by their intersection with the anti-diagonal, 
that is the interval $\langle a, b\rangle = p(\Delta \cap B) $ (as in Proposition~\ref{P:BarofBlock}). Precisely we have:
\begin{lem}\label{L:Blockassociatedtointerval}
 Let  $a<b$ be real numbers. There are unique  blocks $B_{b}^{\langle a, b\rangle}$, $B_{h}^{\langle a, b\rangle}$, $B_{v}^{\langle a, b\rangle}$ and $B_{d}^{\langle a, b\rangle}$ respectively of type $\textbf{bb}^{-}$, $\textbf{hb}$, $\textbf{vb}$ and $\textbf{db}^+$ such that $p(\Delta\cap B_{?}^{\langle a, b\rangle}) = \langle a,b \rangle$. 
\end{lem}
\begin{pf}
By definition~\ref{def:block_MV}, all the blocks except the birth blocks lying entirely in $\Delta^+_{>0}$, that is those of type $\textbf{bb}^{+}$,  are uniquely determined by their intersection with the anti-diagonal (also see figure~\ref{fig:blocks}). In fact the points $a, b$ determine the block of any of these types as in remark~\ref{R:caracofab} and more precisely it determines the boundary lines of the block. To determine if the lines are included in the block or not, we look to whether $a$ or $b$ are inside the interval $\langle a, b\rangle$. For instance, for $B_v^{[a,b]}$ we take the vertical block delimited by the vertical lines $x=-b$ and $x=-a$ and containing them, while $B_v^{]a,b[}$ is the block vertical delimited by the same lines but not containing any of them. 
\end{pf}

\begin{cor}\label{C:pfdimpliesconstructible}
 If $M \in \Obj(\Pe(\Delta^{+}))$ is middle-exact and pointwise finite dimensional, then $\overline{M}$ is weakly constructible. 
 Furthermore, if $M$ is strongly pointwise finite dimensional (definition~\ref{D:spfd}) and midddle-exact, then $\overline{M}$ is constructible. 
\end{cor}
In particular, the restriction of the sheafification functor $\Xi=\overline{(-)}:\Pe(\Delta^{+})\to \Mod(\kk_\R)$ to the full subcategory of pfd modules takes values in the subcategory $ \Mod_{\R c}(\kk_\R)$ of constructible sheaves.
\begin{pf}
 By the decomposition Theorem~\ref{thm:exactdecomp}, the pfd module $M$ is isomorphic to a direct sum of blocks $M\cong  \bigoplus_{B\in \mathbb{B}(M)} \kk^B$. Since $\overline{(-)}$ commutes with direct sum for pfd modules (Proposition~\ref{P:PropertiesofBar}), Proposition~\ref{P:BarofBlock} yields that $\overline{M} \cong \bigoplus_{B\in \mathbb{B}(M)} \overline{\kk^B} $ is a (pointwise finite when $M$ is strongly pfd) direct sum of sheaves of the form $\kk_{I}$ where $I$ is an interval. 
\end{pf}

The level-set persistence to sheaves functor $\Xi$ does not preserve
interleavings in general. However, the trouble is only related to the death or $\textbf{bb}^+$ quadrant. 
More precisely we have the following two lemmas.
\begin{lem}\label{L:Barpreserveentre} Let $M, N\in \Obj(\Pe(\Delta^+))$ be middle exact  pointwise finite dimensional 
and such that their barcodes contains only blocks of type $\textbf{bb}^{-}$, $\textbf{vb}$ and $\textbf{hb}$. 
Then $$\overline{M [\vec{\varepsilon}]} \cong \overline{M}\star K_{\varepsilon}. $$ Furthermore, 
if $M \sim_\varepsilon^{\Delta^+} N$, then 
 $$ \overline{M} \; \sim_\varepsilon \, \overline{N}.$$
\end{lem}
\begin{pf}
 By Theorem~\ref{thm:exactdecomp}, we have  isomorphisms 
 $M \cong\mathop{\bigoplus}\limits_{B\in \mathbb{B}(M)} \kk^B$, 
 $N \cong \mathop{\bigoplus}\limits_{B\in \mathbb{B}(N)} \kk^B$ of persistence modules, 
 such that the blocks $B$ are of types $\textbf{bb}^{-}$, $\textbf{vb}$ and $\textbf{hb}$. Lemma~\ref{L:shiftofBlock} 
 implies that 
 \begin{align*}\overline{M[\vec{\varepsilon}]} &\;\cong \bigoplus_{B\in \mathbb{B}(M)} \overline{\kk^{B -\vec{\varepsilon}}}  \\ 
  \quad \overline{N[\vec{\varepsilon}]} &\;\cong \bigoplus_{B'\in \mathbb{B}(N)} \overline{\kk^{B' -\vec{\varepsilon}} }
  \end{align*}
 where each $\overline{\kk^{B -\vec{\varepsilon}}}$ is of the form 
 $\kk_{I(B, \varepsilon)}$ where $I(B, \varepsilon)$ is an interval $<a,b>=p((B-\vec{\varepsilon})\cap \Delta)$ which is
 \begin{itemize}\item a closed non-empty interval if $B$ is of type $\textbf{bb}^{-}$; 
  \item a semi-open interval closed on the left (resp. closed on the right) if  $B$ is of type $\textbf{hb}$ (resp. $\textbf{vb}$).
 \end{itemize}
Therefore we have: 
\begin{eqnarray*} 
 \text{if $B$ is of type $\textbf{bb}^{-}$, then }  \overline{\kk^{B-\vec{\varepsilon}}} 
 & \cong & \kk_{[a-\varepsilon, b+\varepsilon ]}, \\
 \text{if $B$ is of type \textbf{vquad}, then }  \overline{\kk^{B-\vec{\varepsilon}}} 
 & \cong & \kk_{]a+\varepsilon, b+\varepsilon ]}, \\
 \text{if $B$ is of type \textbf{hquad}, then }  \overline{\kk^{B-\vec{\varepsilon}}} 
 & \cong & \kk_{[a-\varepsilon, b-\varepsilon [}.
\end{eqnarray*}
 Using Proposition~\ref{P:P38Berk18}, we thus get that  in all cases, 
 $$ \overline{\kk^{B-\vec{\varepsilon}}} 
 \; \cong \;  \overline{\kk^B} \star K_\varepsilon $$ and by additivity of the convolution functor we obtain 
 $\overline{M [\vec{\varepsilon}]} \cong \overline{M}\star K_{\varepsilon}$ as claimed.

 \smallskip
 
 The same results holds for  the blocks $B'\in 
 \mathbb{B}(N)$ so that 
 $$\overline{M [\vec{\varepsilon}]} \cong \overline{M}\star K_{\varepsilon}, \quad  \overline{N [\vec{\varepsilon}]}
 \cong \overline{N}\star K_{\varepsilon}.$$
 Note further that, for a  $\textbf{bb}^{-}$ block $B$, the canonical map 
 $\overline{\kk^B}\star K_{\varepsilon} \to \overline{\kk^B}$ (of proposition~\ref{P:propertiesofconvolution}) is identified with the canonical sheaf map 
 $\kk_{[a-\varepsilon, b+\varepsilon ]} \to \kk_{[a, b]}$ as follows from the proof of~\cite[Lemma 3.9]{Berk18}. Since the sheaf map is induced 
 by restriction we obtain from the above equivalences, that the diagram
 \[ \xymatrix{ \overline{\kk^B} \star K_\varepsilon \ar[r]^{\cong}\ar[d] & \kk_{[a-\varepsilon, b+\varepsilon ]} \ar[r]^{\cong} \ar[d] & \overline{\kk^B[\varepsilon]} \ar[d]^{\overline{\tau_{\varepsilon}^{\kk^B}}} 
 \\ \overline{\kk^B} \ar[r]^{\cong} & \kk_{[a, b ]} \ar[r]^{\cong} & \overline{\kk^B}} \] is commutative. Using ~\cite[Proposition 3.10]{Berk18}, 
 the above identification extends to the \textbf{vb} and \textbf{hb} blocks case as well: that is we have, for any block $B$ 
 of type $\textbf{bb}^{-}$, $\textbf{vb}$ and $\textbf{hb}$ a commutative diagram
  \begin{equation}\label{eq:transsenttoconv} \xymatrix{ \overline{\kk^B} \star K_\varepsilon  \ar[r]^{\cong} \ar[d] & \overline{\kk^B[\varepsilon]} \ar[d]^{\overline{\tau_{\varepsilon}^{\kk^B}}} 
 \\ \overline{\kk^B} \ar[r]^{\cong}  & \overline{\kk^B}} \end{equation}
 
 \smallskip
 
 Now let $f: M\to N[\vec{\varepsilon}]$ and $g: N\to M[\vec{\varepsilon}]$ be an $\varepsilon$-interleaving between $M$ and $N$, 
 then applying the functor $\Xi=\overline{(-)}$ to the latter isomorphisms, we obtain  an $\varepsilon$-interleaving in sheaves given by  
 $ \overline{M}\stackrel{\overline{f}}\to \overline{N[\vec{\varepsilon}]}\cong \overline{N}\star K_\varepsilon$ 
 and $\overline{N}\stackrel{\overline{g}}\to \overline{M[\vec{\varepsilon}]}\cong \overline{M}\star K_{\varepsilon}$.
\end{pf}

For \emph{deathblocks} or \emph{birthblocks} of type $\textbf{bb}^{+}$, the sheafification functor does not intertwine shifts 
with convolution in a naive way. 
However we have the following precise result. To state it,  we first 
recall that to a block $B$, we can associate the two real numbers $a,b\in \R$ such that $<a,b> = p(B\cap \Delta)$; 
the convention being that $a=1$ and $b=-1$ when $ p(B\cap \Delta) = \emptyset$. 
\begin{lem}\label{L:Barondeathblocks}Let $B$ be a block of type $\textbf{db}$ or $\textbf{bb}^{+}$. If $B$ is a $\textbf{bb}^{+}$ block,
its dual death block $B^\dagger$ intersects $\Delta$ and we denote $<a^\dagger,b^\dagger> = p(B^\dagger \cap \Delta)$.

Furthermore, for any $\varepsilon \geq 0$, we have that,
 \begin{align}
 &\text{if $B$ is of type \textbf{dquad}, then, } \overline{\kk^{B}[\vec{\varepsilon}] }
 \; \cong \; \left\{ \begin{array}{ll}  \overline{\kk^{B}}\star K_{\varepsilon} & \mbox{ if $\varepsilon< \frac{b-a}{2}$}\\
 0 & \mbox{ if $\varepsilon \geqslant \frac{b-a}{2}$,} \end{array}                                                                                                                          \right .  \label{eq:Barofdb}\\
& \text{if $B$ is of type $\textbf{bb}^{+}$, then, }\overline{\kk^{B}[\vec{\varepsilon}]} 
 \;\cong \;  \left\{ \begin{array}{ll}  0 & \mbox{ if $\varepsilon< \frac{b^\dagger-a^\dagger}{2}$}\\
 \kk_{[\frac{a^\dagger+b^\dagger}{2}, \frac{a^\dagger+b^\dagger}{2}]}\star K_{\varepsilon -\frac{b^\dagger -a^\dagger}{2}} & 
 \mbox{ if $\varepsilon \geqslant \frac{b^\dagger-a^\dagger}{2}$.} \end{array}                                                                                                                          \right .  \label{eq:Barofbb+}
 \end{align}
\end{lem}
\begin{pf}The proof will be similar to the one of lemma~\ref{L:Barpreserveentre}. First note that,  if $B$ is of birthtype $\textbf{bb}^+$,
the supremum of $B^\dagger$ is the infimum of $B$ by definition of the dual block. Therefore, by remark~\ref{R:caracofab}, 
we have that the infimum of the elements of $B$ is the point $(-a^\dagger, b^\dagger)\in \Delta^+$. 
It follows that $B-\vec{\varepsilon}$ remains of type $\textbf{bb}^+$ as long as $\varepsilon < \frac{b^\dagger-a^\dagger}{2}$ 
and it becomes of type $\textbf{bb}^{-}$ when  $\varepsilon \geqslant \frac{b^\dagger-a^\dagger}{2}$. Furthermore, in that latter case, we have that 
$$ p((B-\vec{\varepsilon})\cap \Delta) = \langle b^\dagger  -\varepsilon,  a^\dagger + \varepsilon \rangle.  $$
By lemma~\ref{L:shiftofBlock} and proposition~\ref{P:BarofBlock}, we thus have that {if $B$ is of type $\textbf{bb}^{+}$, then }
 \[
   \overline{\kk^{B-\vec{\varepsilon}}} 
 \, \cong \, \left\{ \begin{array}{ll}  0 & \mbox{ if $\varepsilon< \frac{b^\dagger-a^\dagger}{2}$}\\
 \kk_{[b^\dagger  -\varepsilon,  a^\dagger + \varepsilon]}& \mbox{ if $\varepsilon \geqslant \frac{b^\dagger-a^\dagger}{2}$.} \end{array}                                                                                                                          \right . \]
 Using Proposition~\ref{P:P38Berk18}, we see that for  $\varepsilon \geqslant \frac{b^\dagger-a^\dagger}{2}$, one has $$\kk_{[b^\dagger  -\varepsilon,  a^\dagger + \varepsilon]} 
 \, \cong \,\kk_{[\frac{a^\dagger+b^\dagger}{2}, \frac{a^\dagger+b^\dagger}{2}]}\star K_{\varepsilon -\frac{b^\dagger -a^\dagger}{2}} $$
 which shows the formula~\eqref{eq:Barofbb+}. 
 
 Now, note that if $B$ is of type \textbf{dquad}, then $B-\vec{\varepsilon}$ has a non-empty intersection with $\Delta$ as long as $\varepsilon < \cfrac{b-a}{2}$. And similarly we find, using proposition~\ref{P:BarofBlock} that 
 
 \[
  \overline{\kk^{B-\vec{\varepsilon}}} 
 \, \cong \,\left\{ \begin{array}{ll}  0 & \mbox{ if $\varepsilon \geq \frac{b^\dagger-a^\dagger}{2}$} \\  \kk_{]a+\varepsilon, b-\varepsilon [}  & \mbox{ if $\varepsilon < \frac{b^\dagger-a^\dagger}{2}$}.\end{array}                                                                                                                          \right . \]
 To prove formula~\eqref{eq:Barofdb}, we are left to apply  Proposition~\ref{P:P38Berk18} a last time.
\end{pf}
\begin{remark} The proof of lemma~\ref{L:Barondeathblocks} and proposition~\ref{P:BarofBlock} also shows that the last equivalence in lemma~\ref{L:Barondeathblocks} also reads, for $\varepsilon \geq \frac{b^\dagger-a^\dagger}{2}$, as
 \begin{equation}
  \kk_{[\frac{a^\dagger+b^\dagger}{2}, \frac{a^\dagger+b^\dagger}{2}]}\star K_{\varepsilon -\frac{b^\dagger -a^\dagger}{2}} \; \cong \; \overline{\kk^{B -\overrightarrow{\frac{b^\dagger-a^\dagger}{2}}}} \star  K_{\varepsilon -\frac{b^\dagger -a^\dagger}{2}}.
 \end{equation}

\end{remark}

\section{ Almost isometric equivalence between $\mbox{MV}(\R)^\s$ and $\D^b_{\R c}(\kk_\R)$}
In this section, we explain why the interleaving distance between level set persistence is essentially the same as the derived bottleneck distance between the associated sheaves (in the constructible case).

In order to express this we will relate constructible sheaves by an isometry to a specific type of graded persistence modules, that is those satisfying the following definition.
\begin{defi}\label{D:spfd}
 A middle-exact persistence module $M\in \Pers(\Delta^+)$, is said to be \textbf{strongly pointwise finite dimensional}, if it is pointwise finite dimensional and satisfies the following additional condition : 

$$\text{For every $\alpha \in \Delta$,}~ \varprojlim_{]-x;y[\ni \alpha}M((x,y))  ~ \text{is finite dimensional} $$

A  Mayer-Vietoris system $S=(S_i,\delta^s_i)$ is said to be strongly pointwise finite dimensional if each $S_i$ is strongly pointwise finite dimensional  and only finitely many  $S_i$'s are non-zero.

The \textbf{full subcategory of MV($\R$)} whose objects are \textbf{strongly pointwise finite dimensional MV-systems is denoted by $\text{MV}(\R)^\s$}.
\end{defi}

Our goal now is to build two functors :

$$\xymatrix{
(\overline{-})^{\text{MV}} :&  \text{MV}(\R)^\s \ar[r] &\D^b_{\R c}(\kk_\R)^{\op}  \\
\Psi  :&  \D^b_{\R c}(\kk_\R)^{\op} \ar[r] & \text{MV}(\R)^\s 
}$$

Satisfying for every $F \in \D^b_{\R c }(\kk_\R)$, $(\overline{~\cdot~ })^{\text{MV}} \circ \Psi (F) \simeq F$, in other words $\Psi$ is a pointwise section,
and that, for every $M\in \text{MV}(\R)^\s$, $$d_I^{\text{MV}}\Big(M, \Psi \big( \overline{M}^\MV\big)\Big) = 0.$$
This goal will be achieved by Corollary~\ref{C:MVcircPsiisId} and Corollary~\ref{C:MVdisO}.
\subsection{Construction of the sheafification of MV-systems: the functor $(\overline{~\cdot~ })^{\text{MV}}$}\label{SS:MVfunctor}
We will now apply section~\ref{SS:levelsettopresheaves} to compare the Mayer-Vietoris persistence systems and constructible sheaves. To do so, we first consider the direct sum of the level set persistence to sheaves functor: 

\emph{Let} $\ell:\Mod(\kk_R) \to  \D(\kk_\R)$ \emph{be the localization functor} sending the category of complexes of sheaves over $\R$ to its derived category $\D(\kk_\R)$.
\begin{defi}\label{D:MVtoSheaves} The \emph{sheafification of MV-systems functor}: $\overline{(-)}^{MV}: \text{M-V}(\R) \to  \D(\kk_\R)^{\op}$ is the functor given, on
objects $S=(S_i, \delta_i^S)_{i\in \Z, s\in \R^2_{>0}}\in \Obj(\text{M-V}(\R))$, by 
$$ \overline{S}^{MV} := \ell\Big(\bigoplus_{i\in \Z} \overline{S_i}[-i]\Big)$$ and, on morphisms $(S_i\stackrel{\varphi_i}\to T_i)_{i\in \Z}$, by 
$$\overline{(\varphi_i)_{i\in \Z}} := \ell\left(\bigoplus \overline{\varphi_i}\right). $$ 
\end{defi}
That this is a functor is a direct consequence of section~\ref{SS:levelsettopresheaves}.

\begin{lem}\label{L:MVofsfrestricts}
 If $S$ is a strongly pointwise finite dimensional Mayer-Vietoris system, then $\overline{S}^{MV}$ is a constructible sheaf. In particular, we have a commutative diagram of functors:
 $$\xymatrix{ \text{M-V}(\R) \ar[r]^{\quad\overline{(-)}^{MV}} & \D(\kk_\R)^{\op} \\ 
 \text{M-V}(\R)^{\s} \ar@{^{(}->}[u]\ar[r]^{\quad\overline{(-)}^{MV}} & \D^b_{\R c}(\kk_\R)^{\op} \ar@{^{(}->}[u].}  $$
\end{lem}
\begin{pf}
 We can apply Theorem~\ref{thm:decom_MV} together with proposition~\ref{P:BarofBlock} in a way similar to the proof of corollary~\ref{C:pfdimpliesconstructible}.
\end{pf}
The same argument shows that if $S$ is pfd (but not necessarily strongly), then $\overline{S}^{MV}$ is weakly constructible.

\begin{prop}\label{P:PropertiesofMVSheafification}
The sheafification of MV-systems functor $\overline{(-)}^{MV}: \text{M-V}(\R)^{\s} \to  \D^b_{\R c}(\kk_\R)$ satisfies the following properties : 
\begin{enumerate}
    \item it commutes with degree shifting operator : for all Mayer-Vietoris system $M$, one has $\overline{M[n]}^{MV} \cong \overline{M}^{MV}[n]$. 
    \item For  a block $B$ of type \textbf{bb$^-$}, \textbf{hb}, \textbf{vb}, \textbf{db}$^{+}$, $j\in \Z$ and $\varepsilon \geq 0$, we have : $$(\overline{S_j^B [\varepsilon]})^{\MV} \simeq \kk_{I(B)}[-j] \star K_\varepsilon $$ 
    where, still denoting  $\langle a, b\rangle=p(B\cap \Delta^+)$, $I(B)$ is the interval given by \begin{eqnarray*} I(B)= [a,b] \mbox{ if $B$ is of type \textbf{bb}$^{-}$,} & &  I(B)= [a,b[ 
    \mbox{ if $B$ is of type \textbf{hb},} \\
      I(B)= ]a,b] \mbox{ if $B$ is of type \textbf{vb},} & &  I(B)= ]a,b[ \mbox{ if $B$ is of type \textbf{db}$^{+}$.}                                                                                                                                                                                      \end{eqnarray*}

    \item If $M \sim_{\varepsilon}^{\Delta^+} N$, then $\overline{M}^{MV} \sim_{\varepsilon} \overline{N}^{MV}$.
    \item If $\overline{M}^{MV}$ is isomorphic to  $\overline{N}^{MV}$ (in the derived category), then $d_I^{MV}(M, N) =0$.
 \end{enumerate}
\end{prop}
\begin{pf}
 Note that assertion 1 is immediate from the definition since we put each $\overline{S_i}$ precisely in degree $i$.

\smallskip 

2 and 3. First assume $B$ is of type  \textbf{bb$^-$}, \textbf{hb} or \textbf{vb}. Then  definition~\ref{D:blocksmodulesforMV} implies that $S_j^B \cong \kk_B[-j]$. Since  $\overline{(\cdot)}^{MV}$ commutes with direct 
sum and shifts, Lemma~\ref{L:Barpreserveentre} implies $(\overline{S_j^B [\varepsilon]})^{\MV} \simeq \kk_B[-j] \star K_\varepsilon $ and further, 
for $\varepsilon \leq \varepsilon'$, 
this isomorphism sends the canonical structure maps $S_j^B [\varepsilon']\to S_j^B [\varepsilon]$ onto the canonical map  
$\kk_B[-j] \star K_\varepsilon \to \kk_B[-j] \star K_{\varepsilon '}$ (see diagram~\eqref{eq:transsenttoconv}).

It remains to prove the same result in the case of a block of type \textbf{db}. Then definition~\ref{D:blocksmodulesforMV} says that as a graded persistent module, one has $$S_j^B \cong \kk_B[-j] \oplus \kk_{B^\dagger}[-j-1]$$ and therefore  $$\overline{S_j^B[\vec{\varepsilon}]}^{MV} \cong \overline{\kk_B[\varepsilon]}[-j] \oplus \overline{\kk_{B^\dagger}[\varepsilon]}[-j-1].$$ Denote $\langle a, b\rangle =p(B \cap \Delta)$ as before Proposition~\ref{P:BarofBlock}. Following the notation of 
Lemma~\ref{L:Barondeathblocks} we thus have that for the dual block $B^\dagger$ of type \textbf{bb}$^{+}$, one has that $a^\dagger=a$, $b^\dagger=b$ by definition. Then, Lemma~\ref{L:Barondeathblocks}, the commutation of convolution wih shifts and Proposition~\ref{P:BarofBlock} imply that 
\begin{multline}\label{eq:MVonSofDB}
\overline{\kk_B[\varepsilon]}[-j] \oplus \overline{\kk_{B^\dagger}[\varepsilon]}[-j-1] \;
  \cong \;\left\{ \begin{array}{ll}  \kk_{]a,b[}\star K_{\varepsilon}[-j] & \mbox{ if $\varepsilon< \frac{b-a}{2}$}\\
 \kk_{[\frac{a+b}{2}, \frac{a+b}{2}]}\star K_{\varepsilon -\frac{b -a}{2}} [-j-1] & \mbox{ if $\varepsilon \geqslant \frac{b-a}{2}$.} \end{array}   \right .
\end{multline}
This formula ~\eqref{eq:MVonSofDB} is precisely the formula for $\kk_{]a,b[}\star K_{\varepsilon}$ according to Proposition~\ref{P:P38Berk18}. We obtain 
a commutative diagram similar to~\eqref{eq:transsenttoconv} in the same way as in Lemma~\ref{L:Barpreserveentre}.
This concludes the proof of claim 2.
Assertion 3 follows immediately of assertion 2 and the fact that the canonical translation maps of persistent modules are sent to the canonical maps 
$\kk_B[-j] \star K_\varepsilon \to \kk_B[-j] \star K_{\varepsilon '}$.

\smallskip

4. Assume $\overline{M}^{MV} \cong \overline{N}^{MV}$. By Theorem~\ref{thm:decom_MV}, we can decompose 
$$M\cong \mathop{\bigoplus}\limits_{j\in \Z} \left(\mathop{\bigoplus}\limits_{B_M \in \mathbb{B}_j(M)}S_j^{B_M}[-j]\right) \mbox{ and }
N \cong  \mathop{\bigoplus}\limits_{j\in \Z} \left(\mathop{\bigoplus}\limits_{B_N \in \mathbb{B}_j(N)}S_j^{B_N}[-j]\right)$$ into  Mayer-Vietoris blocks.
Since $\overline{(-)}^{MV}$ commutes with direct 
sum and shifts (by property 1),  we have isomorphisms
\begin{eqnarray*}
 \overline{ \bigoplus_{j\in \Z} \bigoplus_{B_M\in \mathbb{B}_j(M)}S_j^{B_M}[-j]}^{MV} &\cong&
 \overline{\bigoplus_{j\in \Z} \bigoplus_{B_N\in \mathbb{B}_j(N)}S_j^{B_N}[-j] }^{MV}  \\ 
 \bigoplus_{j\in \Z} \bigoplus_{B_M\in \mathbb{B}_j(M)} \overline{S_j^{B_M}}^{MV}[-j] &\cong& 
 \bigoplus_{j\in \Z} \bigoplus_{B_N\in \mathbb{B}_j(N)} \overline{S_j^{B_N} }^{MV}[-j].
\end{eqnarray*}
For any vertical, horizontal or \textbf{bb}$^-$ type block $B$, Proposition~\ref{P:BarofBlock} tells us that  
$\overline{S_j^{B}}^{MV} \cong \kk_{I(B)}$ where $I(B)$ is a non-empty interval (uniquely determined by $p(B\cap \Delta)$). 
If $B$ is of type \textbf{db}$^{+}$, then 
$$
 \overline{S_j^{B}} \; \cong \;  \kk_{I(B)} \oplus \kk_{I(B^\dag)}[-1]
$$
according to definition~\ref{D:blocksmodulesforMV} and~\ref{D:MVtoSheaves}.  
Therefore, we have an isomomorphism 
\begin{multline}
 \label{eq:decofBarMV}  
 \bigoplus_{j\in \Z} \left(\Big(\bigoplus_{B_M\in \mathbb{B}_j(M)\setminus \mathbb{B}^{\textbf{dq}}_j(M)} \kk_{I(B_M)}[-j] \Big) \oplus 
 \Big(\bigoplus_{B_M\in \mathbb{B}^{\textbf{dq}}_j(M) } \big(\kk_{I(B_M)}[-j] \oplus \kk_{I(B_M^\dag)}[-j-1]\Big)\right)\\
 \cong \; 
 \bigoplus_{j\in \Z} \left(\Big(\bigoplus_{B_N\in \mathbb{B}_j(N) \setminus \mathbb{B}^{\textbf{dq}}_j(N)} \kk_{I(B_N)}[-j]\Big)\oplus 
\Big( \bigoplus_{B_N\in \mathbb{B}^{\textbf{dq}}_j(N) } \big(\kk_{I(B_N)}[-j] \oplus \kk_{I(B_N^\dag)}[-j-1]\Big)\right)   
\end{multline}
of constructible sheaves. Here $\mathbb{B}^{\textbf{dq}}_j(M)$, $\mathbb{B}^{\textbf{dq}}_j(N) $ are the subsets of those bars that are of type \textbf{db}$^{+}$ in the respective decompositions of $M$ and $N$. 

By unicity of the decomposition in Theorem~\ref{T:KSdecomposition}, we obtain  degreewise bijections between 
the set of associated graded barcodes $\{ I(B_M), \, B_M \in  \mathbb{B}_j(M)\}$ and $\{ I(B_N), \, B_N \in  \mathbb{B}_j(N)\}$
and therefore  bijections $\sigma_j: \mathbb{B}_j(M) \cong \mathbb{B}_j(N)$ with the property that  for any $B_M\in \mathbb{B}_j(M)$, 
$\sigma_j(B_M)$ is a block of the same type as $B_M$ and which is equal to $B_M$ except maybe on the boundary.

\begin{lem} Let $\mathcal{B}, \mathcal{B}'$ be sets of M-V blocks of types \textbf{db}, \textbf{vb}, \textbf{db} and \textbf{bb}$^-$. 
If there is a bijection $\sigma: \mathcal{B} \to \mathcal{B}'$  such that for any $B\in \mathcal{B}$, $\sigma(B)$ is equal to $B$ except maybe on the boundary, then 
$$d_I^{MV}\left(\bigoplus_{B\in \mathcal{B}} S^B_{j}, \bigoplus_{B'\in \mathcal{B}'} S^{B'}_{j}\right) =0.$$
\end{lem}
\begin{pfofLemma} It is enough to check that,
 if $B$ and $B'$ are two blocks of the same type which differs only on their boundary, then $B$ and $B'$ are $\varepsilon$-interleaved for any $\varepsilon >0$. This property follows from Lemma~\ref{L:shiftofBlock} and an immediate application of the definition of the blocks of each type. 
 Then  the direct sum of those interleavings relating each $\mathcal{B}$ to $\sigma(\mathcal{B})$  gives a $\varepsilon$-interleaving in between $\bigoplus_{B\in \mathcal{B}} S^B_{j_B}$ and  $\bigoplus_{B'\in \mathcal{B}'} S^{B'}_{j_{B'}}$ for every $\varepsilon >0$; the lemma follows.
\end{pfofLemma}
The  claimed property 3 follows from the lemma since we have proved just above that we can find such a 
permutation relating $\mathbb{B}_j(M)$, $\mathbb{B}_j(N)$ for each degree $j$.
\end{pf}

Let $f:X\to \R$ be a continuous map. Then we have the derived functors of the direct image: $\Rr^i f_* \kk_X \in \Mod(\kk)$, see~\cite{Kash90, Iversen} which are the cohomology groups of the derived functor $\Rr f_* \kk_X \in \D(\kk_X)$. Note that this is just a special case of derived direct image, defined for any continuous map $\phi: X\to Y$, which   is a functor $\Rr \phi_*: \D(\kk_X) \to \D(\kk_Y)$. In particular, the\begin{equation}\label{eq:defofdirectimageasfunctor} \textit{assignment } f\mapsto  \mathop{\bigoplus}\limits_{i\in \N} \Rr^i f_*(\kk_X)[-i] \textit{ defines a functor } \Rr (-) \kk_{(-)}: \mathrm{Top}_{| \R} \to \D(\kk)^{\op}.\end{equation} 
A morphism $\phi: (X,f) \to (Y,g)$ is mapped by this functor to the linear map $\bigoplus \Rr^i\phi_*: \Rr^i g_* \kk_Y \to \Rr^i f_* \kk_X$ and the fact that this defines a functor is an immediate consequence of the
composition formula $\Rr (\psi\circ \kappa)_* \; \cong \; \Rr \psi_* \circ \Rr \kappa_* $ see~\cite{Kash90, Iversen}.
\begin{prop}\label{P:MVmapstoRf} Assume $X$ is locally contractible. Then there is an natural isomorphism $$\overline{M^f}^{MV} \; \cong \; \bigoplus_{i\in \N} \Rr^if_* \kk_X [-i].$$
\end{prop}
\begin{pf}
 By example~\ref{Ex:MVfromFct} and definition~\ref{D:MVtoSheaves}, we have 
 $ \overline{M^f}^{MV} \; \cong  \; \bigoplus_{i\in \N} \ell (\overline{H_i(f^{-1}(-))}) [-i] $. 
 Now, from definition~\ref{D:defofBar}, we have that $\overline{H_i(f^{-1}(-))}$ is the sheafification of the presheaf 
 \[
 \Ouv(\R)\ni U \mapsto  \varprojlim_{]-x,y[\subset U} \Hom_{\kk}\left(H_i(f^{-1}(]-x,y[)) ,\kk\right) \; \cong\; H^i(f^{-1}(U))
 \]
since $\kk$ is a field (and therefore the cohomology of the dual of a chain complex is the dual of the homology) and every open in $\R$ is a disjoint union of intervals. It is well-known that for $f:X\to Y$ and any sheaf $F$,  $\Rr^i f_*(F)$ is the sheaf associated to the presheaf $\Ouv(Y)\ni U \mapsto  H^i(f^{-1}(U),F)$ (see \cite[Proposition 5.11]{Iversen} for instance). Furthermore, when $X$ is locally contractible, one has an isomorphism of presheaves $$V\mapsto H^i(V,\kk_X) \, \cong \,H^i(V, \kk_V) \, \cong\, H^i(V)$$ where the first isomorphism is for the sheaf cohomology with value in a constant sheaf and its restriction ${\kk_X}_{|_V} \cong \kk_V$ to an open subset, and the last isomorphism is the usual identification of sheaf cohomology with value on a constant sheaf with singular cohomology for locally contractible spaces.
\end{pf}

\subsection{The functor $\Psi$ from constructible sheaves to Mayer-Vietoris systems}\label{SS:Psi}
We now turn to the construction of a section of the sheafification of MV-systems. We have the following intrinsic definition. 
\begin{defi}\label{def:Psi}Given $F^\bullet \in \D \big(\Mod(\kk_{\mathbb{R}})\big)$ and $i \in \Z$, we define  $\Psi(F^\bullet)_i$ to be the object of $\Pers(\Delta^+)$ given, for $(x,y)\in \Delta^+$ by :
\begin{equation}\label{eq:defPsiHypercohom}
 \Psi(F^\bullet)_i \; = \; \Hom_{\Mod(\kk)}\left ( \mathbb{H}^i\Big(F^\bullet_{|\, ]-x,y[}\Big), \kk\right)
\end{equation}
where $\mathbb{H}^*(-)$ is the hypercohomology of  complexes of sheaves and $F^\bullet_{|\, ]-x,y[}$ is the restriction of $F^\bullet$ to the open $]-x, y[$.  
\end{defi} 
The degree $i$ hypercohomology of the restriction $F^\bullet_{|U}$ of $F^\bullet$ to an \emph{open} set $U$ is isomorphic to the hypercohomology $\mathbb{H}^i(U,F^\bullet)$ of $F^\bullet$ on the open $U$ (see~\cite{Kash90} for instance) so that the latter formula can be used in~\eqref{eq:defPsiHypercohom}.

There is a slightly less homological algebraic involved formula for constructible sheaves, see Lemma~\ref{L:DefadnFunctorialityPsi} below.  Since
those are our case of interest, the reader can take  formula~\eqref{eq:Psifocosntructible} as the definition of $\Psi_i$ for the rest of the paper. 

\smallskip

That $\Psi(F^\bullet)_i$ is a persistent module follows from Lemma~\ref{L:Defiota}. Since hypercohomology and duality are functors, then it is immediate that for all $i\in \Z$,
\begin{equation}
 \Psi(-)_i: \D \big(\Mod(\kk_{\mathbb{R}})\big)^{\op} \longrightarrow \Pers(\Delta^+) \mbox{ is a functor.}
\end{equation}
 Note that the hypercohomology (see~\cite{DanilovHyperhomology, Tohoku} for standard references) of a complex of sheaves $G^\bullet$ on a space $X$ is obtained by replacing $G^\bullet$ by a quasi-isomorphic \emph{injective complex} of sheaves\footnote{that is a fibrant resolution in the model category of sheaves} $I^\bullet$ (which in the case where $F^\bullet$ is bounded on the left is the same as a quasi-isomorphic chain complex whose terms are injective sheaves) and then taking the cohomology of the section of this  complex $I^\bullet$: 
  $$\mathbb{H}^i(G^\bullet)\; =\; H^i( \Gamma(I^\bullet)) .$$

  \smallskip 
  
For constructible sheaves $\bigoplus \kk_{I}[n_I]$, one can simply compute the hypercohomology by  computing the derived sections of the sheaf as follows immediately from the following lemma (since the homology is equal to the sheaf in that case).  
\begin{lem}\label{L:DefadnFunctorialityPsi}
 Let $F$ be in $\D^b_{\R c}(\kk_\R)$. Then for any $i\in \Z$, and $(x,y)\in \Delta^+$, one has an   isomorphism of persistent modules over $\Delta^+$:
 \begin{equation}\label{eq:Psifocosntructible}\Psi(F)_i(x,y)\; \cong\; \bigoplus_{k+l = i} \Hom_{\Mod(\kk)}\left ( \Rr^k\Gamma \left ( ]-x,y[ ,  H^l(F)\right ), \kk \right )  \end{equation}
 where $\Rr^k\Gamma \left ( ]-x,y[ , -\right )$ is the $k$-th right derived functor of the functor of sections on $]-x,y[$ and $H^*(F)$ is the graded sheaf given by the homology of the underlying complex of $F$. 
\end{lem}
Note that since $F$ is assumed to be constructible, \emph{there are only finitely many} pairs $(k,l)$ such that the right-hand-side vector space is non zero.

  \begin{pf}
  The reader who knows the spectral sequences associated to hypercohomology can immediately deduce the result of the lemma by noticing 
  that the assumption on $F$ implies its degeneracy at the $E_2$-page which is exactly the right hand side of~\eqref{eq:Psifocosntructible}.
   
   \smallskip
   
   Alternatively, let $F$ be any complex of sheaves on a space $X$. 
   Denote $\Gamma(-,F): \Ouv(X)^{\op}\to \Mod(\kk)$ the functor sending an open $U$ to the sections $F(U)$ of $F$ over $U$. Denote 
   $\Rr^*\Gamma(-,F): \Ouv(X)^{\op}\to \D(\Mod(\kk))$ its derived functor, which, by definition is given by $\Gamma(-, I^\bullet)$ 
   where $I^\bullet$ is an injective complex of sheaves quasi-isomorphic to $F$.  Note also that 
   $\Rr^* \Gamma(U, F) \cong \Rr^*(F_{|U})$ and that for any sheaf $F$, one has $\Rr^k(U,F)= H^k(\Rr^*\Gamma(U,F)$, 
   see~\cite{Kash90} or another classical textbook.
   
   Then according to definition~\ref{def:Psi} $\Psi_i(F)$ is the persistent object given by the composite functor
  \begin{equation}\label{eq:Psiascomposite}\xymatrix{\Delta^+ \ar[r]^{\iota\qquad} &\Ouv(\R) \ar[rr]^{\Rr^*\Gamma(-,F)\quad } && \D(\Mod(\kk))^{\op} \ar[r]^{\;\;H^i(-)} & \Mod(\kk)^{\op} \ar[rr]^{\Hom_{\kk}( - , \kk)} && \Mod(\kk).}\end{equation}
  Now we asume $F \in \D^b_{\R c}(\kk_\R)$. By Theorem~\ref{T:KSstructure}, we have an isomorphism of complexes of sheaves 
  $F \cong \bigoplus_{j\in \Z} H^j(F)[-j]$. Therefore we can replace $F$ by its homology in~\eqref{eq:Psiascomposite}. Then,  
  we can take $I^\bullet$ to be the direct sum of injective resolutions of each $ H^j(F)[-j]$. 
  The lemma follows  thanks to the fact that only finitely many $k$ and $l$ in~\eqref{eq:Psifocosntructible} gives non-zero terms as noted above and therefore 
  the functors in~\eqref{eq:Psiascomposite} commutes with the (finite) direct sum.
  \end{pf}
\begin{remark}\label{rk:Psinotfaithful} The functor $\Psi$ is not faithful. Indeed for any $a<b<c$, one has an non-split exact sequence of sheaves
\[0\to \kk_{[a,b[}\to \kk_{[a,c[} \to \kk_{[b,c[}\to 0 \] which gives a non zero homomorphism $\kk_{[b,c[}\to \kk_{[a,b[}[1]$ in  $\D^b_{\R c}(\kk_\R)$.
However there are no non-zero Mayer-Vietoris systems homomorphism in between $\Psi(\kk_{[b,c[})$ and $\Psi(\kk_{[a,b[}[1])=\Psi(\kk_{[a,b[})[1]$ 
as follows from Proposition~\ref{P:PsionIntervals} below since there are non non-zero  homomorphims in between M-V systems associated to 
horizontal blocks in \emph{different} degrees. 

\smallskip

 Note also that the isomorphism of lemma~\ref{L:DefadnFunctorialityPsi} is \emph{not} natural in $F$ for similar reasons. For instance, the right hand side of~\eqref{eq:Psifocosntructible} maps the non zero morphism  $\kk_{[b,c[} \to \kk_{]a,b[} [1]$ (induced by the short exact sequence $0\to \kk_{]a,b[} \to\kk_{]a,c[}\to \kk_{[b,c[} \to 0$) to $0$ but $\Psi$ does not. 
\end{remark}

\begin{remark}
 The functor $\Psi$ is thus essentially defined as  the dual of the derived section of $F$ and not just as the dual of the homology sheaf of 
 $F$ which could have been a more naive approach. The main reason is that the latter will not carry a Mayer-Vietoris structure; in other words, 
 it will forget too much of the structure of the constructible sheaf. However, the derived construction carries such a structure in a natural way 
 as we will now see.
\end{remark}

\begin{prop}\label{P:psiofCOnstrisMV}
The family $(\Psi(F)_i)_{i\in \Z}$ carries a natural structure of a Mayer-Vietoris system. In addition, if $F\in \D^b_{\R c}(\kk_\R)$,
then  it is strongly pointwise finite dimensional (Definition~\ref{D:spfd}). 
\end{prop}

\begin{pf} We have already seen  that $\Psi_i(F)$ is a persistence module over $\Delta^+$ as an immediate consequence of lemma~\ref{L:Defiota}.
For $s = (s_1,s_2)\in \R^2_{>0}$ and $i\in\Z$, we have to build the connection morphism $\delta_i^2$. Let $I^\bullet\in C^b(\kk_\R)$
an injective resolution of $F$ in the category of  sheaves. Consider $(x,y)\in \Delta^+$, 
then we have the Mayer-Vietoris sequence associated to the cover $]-x - s_1 , y[ \cup ]-x , y + s_2[ $ of $]-x - s_1 , y + s_2[$ 
which is the short exact sequence of complexes of sheaves 

$$0 \longrightarrow \Gamma(]-x - s_1 , y + s_2[, I^\bullet) \longrightarrow \Gamma(]-x - s_1 , y[, I^\bullet) \oplus  \Gamma(]-x , y+s_2[, I^\bullet) 
\longrightarrow \Gamma(]-x,y[, I^\bullet) \longrightarrow 0 .$$
Let us write $H^i(U, I^\bullet)$ for the $i$-th cohomology groups $\Rr^i(U, I^\bullet)$. 
Passing to cohomology, we thus obtain a long exact sequence (see~\cite{Kash90}) 
\begin{multline}\label{eq:lessheavesforPi} 
\dots \to     H^i(]-x-s_1, y+s_2[, I^\bullet) \to 
 H^i(]-x-s_1, y[, I^\bullet) \oplus  H^i(]-x, y+s_2[, I^\bullet) \\ \to 
  H^i(]-x, y[, I^\bullet) \stackrel{\delta}\to H^{i+1}(]-x-s_1, y+s_2[, I^\bullet) \to \dots.
\end{multline}
Since by definition of sheaf cohomology, one has, $H^i(]-x, y[, I^\bullet) \cong \Rr^i\Gamma \left ( (]x,y[ , F\right )$,
the linear dual of the maps $\delta$ given by the exact sequence~\eqref{eq:lessheavesforPi} yields linear maps 
$\delta_i^s:\psi_i(F)((x,y)) \to \psi(F)[s](x,y)$ for all $(x,y) \in \Delta^+$. The exactness of~\eqref{eq:lessheavesforPi} and Lemma~\ref{L:Defiota} 
also implies that the collection $(\psi_i(F), \delta_i^s)_{i,s}$ is a Mayer-Vietoris system over $\R$.

\smallskip 

When $F$ is constructible, its cohomology groups are finite dimensional in each degree, and there are only finitely many of them. 
Therefore $\Psi(F)$ is pointwise finite dimensional. Now the proof that $\Psi(F)$ is strongly finite dimensional is an argument similar 
to the proof of property 4 in Proposition~\ref{P:PropertiesofBar}. Alternatively, one can simply use the structure theorem~\ref{T:KSstructure}
and proposition~\ref{P:PsionIntervals} below to conlude directly since  strongly pointwise finite dimensional modules are stable under 
locally finite direct sums.\end{pf}

\begin{prop}\label{P:PropertiesofPsi}
The rule $F\mapsto \Psi(F):= (\Psi_i(F), \delta_i^s)_{i,s}$ defines  functors $\Psi: \D^b_{\R c}(\kk_\R)^{\op}\to \MV(\R)^\s$,  $\Psi: \D(\kk_\R)^{\op}\to \MV(\R)$ fitting in a commutative diagram: 
\[\xymatrix{ \D(\kk_\R)^{\op}\ar[r]^{\Psi} & \MV(\R)\\ 
\D^b_{\R c}(\kk_\R)^{\op}\ar[r]^{\Psi} \ar@{^{(}->}[u]& \MV(\R)^\s  \ar@{^{(}->}[u].} \]
Furthermore, these functors are additive and commutes with shifts associated to the canonical triangulated structure of the derived category.
\end{prop}
\begin{pf} Since the definition of $\Psi_i$ is functorial  and the connecting morphism in Mayer-Vietoris long exact sequences is also functorial, we obtain  that $\Psi$ is indeed a functor. Proposition~\ref{P:psiofCOnstrisMV} gives the fact that $\Psi$ sends the subcategory of constructible sheaves to the one of strongly pointwise finite dimensional systems.  The last assertion follows from the fact that hyperchohomology commutes with direct sums and shifts. 
\end{pf}
\begin{ex}\label{ex:ValueofPsiondirectsum}
 Let $F =\bigoplus_{I} \kk_I[n_I]$ be constructible (derived) sheaf over $\R$. Then by Proposition~\ref{P:PropertiesofMVSheafification} and Lemma~\ref{L:DefadnFunctorialityPsi} we obtain, for any $(x,y)\in \Delta^+$, the simple formula 
 $$\Psi(F)_i (x,y) ) \, \cong \, \bigoplus_{I} \bigoplus_{k} \Rr^k\Gamma\left(]-,x,y[, \kk_I\right)[n_i+k]$$ for $\Psi(F)$ (one can also note that the only values of $k$ for which we have a non zero term are $0$ and $1$ from Proposition~\ref{P:BarofBlock}).
\end{ex}

Recall definition~\ref{D:blocksmodulesforMV} of the canonical MV-systems $S_i^B$ associated to blocks as well as Lemma~\ref{L:Blockassociatedtointerval}. 
The following is the analogue of proposition~\ref{P:BarofBlock}, that is, it describes the action of $\psi$ on the building blocks of a constructible sheaf. Together with example~\ref{ex:ValueofPsiondirectsum}, it allows to compute the value of $\Psi$ explicitly. 
\begin{prop}\label{P:PsionIntervals}  Let $I=\langle a, b\rangle$ be an interval in $\R$. 
\begin{enumerate}
    \item If $I$ is open, then  $\Psi ({\kk_{]a, b[}}[-i]) \cong 
S_i^{B_{d}^{[ a, b ]}}$.
    \item If $I=]a,b]$, then  $\Psi ({ \kk_{] a, b]}}[-i]) \cong 
S_i^{B_{v}^{[ a, b[}}$.
    \item If $I=[a, b[$, then  $\Psi ({ \kk_{[ a, b[}}[-i]) \cong 
S_i^{B_{h}^{]a, b]}}$. 
    \item If $I$ is compact, then  $\Psi ({ \kk_{[a, b]}}[-i]) \cong 
S_i^{B_{b}^{]a, b[}}$. 
\end{enumerate}
Here all the isomorphisms are isomorphisms of Mayer-Vietoris systems and the blocks $S_{?}^{\langle a, b\rangle}$ are given by lemma~\ref{L:Blockassociatedtointerval}. 
\end{prop}
\begin{remark}
 Note that when applying $\Psi$ on an interval, the closed boundary becomes an open boundary lines in the associated block of the image and the open ones  become closed.
\end{remark}

As will be made clear by the proof, the claim 1  relies heavily on the fact that we have taken a derived functor approach for the definiton of $\Psi$. 

\begin{pf}Let us first prove the open interval case.
In view of the proof of proposition~\ref{P:psiofCOnstrisMV}, using compatibility with shifts and direct sums, we only need to compute the cohomology groups of $\Rr^k\Gamma \left ( ]-x,y[ , \kk_{I}\right )$ which by definition (see~\cite{Kash90}) is isomorphic to 
$Ext^k_{\kk_\R}\left( \kk_{]-x,y[}, \kk_{I }\right) $. 
By Proposition 3.13 and 3.14 in~\cite{Berk18}, we have that it is always $0$ for $k>1$. Furthermore, the only case for which it is non-zero for $k=1$ is when $I$ is an open whose closure is included in $]-x,y[$. In that latter case (which means, if $I=]a,b[$, that $[a,b]\subset ]-x,y[$ i.e. $x>-a$ and $y>b$) we  then have $Ext^1_{\kk_\R}\left( \kk_{]-x,y[}, \kk_{I }\right) \cong \kk $.  Therefore, by functoriality of the $Ext^1_{\kk_\R}( -, \kk_I)$ functor in its left variable, it follows that the persistence module associated to  
$Ext^1_{\kk_\R}( -, \kk_I)$ in $\Psi(\kk_I)$ is either $0$ if $I$ is not open or, if $I$ is open, is precisely the block module in degree $1$ which is supported on the type \textbf{bb}$^+$ block whose infimum is $(-a,b)$ and contains none of its boundary lines. Here, by block module we refer to Definition~\ref{Def:blockmodule}. Therefore, by definition of duality~\ref{D:dualblocks}, for an open $I=]a,b[$ , the contribution of $Ext^1_{\kk_\R}( -, \kk_I)$ in $\Psi(\kk_I)$ is precisely 
$\kk^{(B_d^{[a,b]})\dag}[-1]$ in degree 1 supported on the type $\textbf{bb}^{+}$ block dual to the deathblock $B_d^{[a,b]}$. 

It remains to compute the image of the $Ext^0_{\kk_\R}\left( \kk_{]-x,y[}, \kk_{I }\right) $. By Proposition 3.13 and 3.1 in~\cite{Berk18}, we find that if $I$ is open, 
$$Ext^0_{\kk_\R}\left( \kk_{]-x,y[}, \kk_{I }\right) \cong \left \{ \begin{array}{ll}
\kk & \mbox{if } ]-x, y[ \subset I \\ 0 & \mbox{else. }
\end{array}\right .$$
For $I=]a,b[$, the condition $]-x,y[ \subset ]a,b[$ can be rewritten as $x\leq -a$ and $y\leq b$.
Using functoriality of $Ext$ again, we thus find that, when $I$ is open, the persistence module associated to $\kk_I$ is the block module $ \kk^{B_d^{[a,b]}}$ concentrated in degree $0$ and supported on the type $\textbf{db}$ block $B_d^I$. Combining the degree $0$ and $1$ part, 
the functoriality of the the Mayer-Vietoris long exact sequence~\eqref{eq:lessheavesforPi} then shows that $\Psi(\kk_I) $ is precisely the MV-block module $S_0^{B_d^{[a,b]}}$ as in Definition~\ref{def:block_MV}.

\smallskip

Now for the three other types of intervals, the computation is easier since we only have to consider $Ext^0_{\kk_\R}\left( \kk_{]-x,y[}, \kk_{I }\right)$ in the computation of $\Psi(\kk_I)$ (all other degrees are $0$ by the $Ext$ computations of~\cite{Berk18}). Arguing as for the open interval case, using Proposition 3.13 and 3.1 in~\cite{Berk18}, we obtain  that the persistence modules $Ext^0_{\kk_\R}\left( \kk_{-}, \kk_{I }\right)$ are respectively the block modules $S_0^{B_v^{[a,b[}}$, $S_0^{B_h^{]a,b]}}$ and 
$S_0^{B_b^{]a,b[}}$ when $I$ is of the type $]a, b]$, $[a, b[$ or $[a,b]$. 
\end{pf}

\begin{prop}\label{P:Psipreservesinterl}Let $F \in \D^b_{\R c}(\kk_\R)$. For any $\varepsilon \geq 0$, there is an isomorphism of graded persistence modules 
$ \Psi(F \star K_{\varepsilon}) \; \cong \; \Psi(F)[\vec{\varepsilon}]$.
\end{prop}
\begin{pf}
 Using theorem~\ref{T:KSdecomposition}, we have that $F\cong \bigoplus_{\alpha \in A} \kk_{I_\alpha} $. By compatibility of convolution with direct sums and shifts, it is thus enough to prove the result for $\kk_I$ for an interval $I$. 
 
Let us start with the case where $I= [a,b]$ is compact.    Then by Proposition~\ref{P:P38Berk18}, we obtain 
\begin{equation}\label{eq:Psiofcpctshift}
 \Psi(\kk_{[a,b]}[-i] \star K_{\varepsilon}) \,\cong\,  \Psi(\kk_{[a-\varepsilon,b+\varepsilon]}[-i])
 \, \cong \, S_i^{B_{b}^{\langle a-\varepsilon, b +\varepsilon\rangle}}
\end{equation}
where the last isomorphism is given by Proposition~\ref{P:PsionIntervals}. Note that, by definition, the block $B_b$ is of type \textbf{bb}$^{-}$ (see Lemma~\ref{L:Blockassociatedtointerval} and definition~\ref{D:blocksmodulesforMV}). Therefore as a persistent module over $\Delta_+$, we have 
$S_i^{B_{b}^{\langle a-\varepsilon, b +\varepsilon\rangle}} \,\cong\, \kk^{B_{b}^{\langle a-\varepsilon, b+\varepsilon \rangle}} [-i]$. By lemma~\ref{L:shiftofBlock} and remark~\ref{R:caracofab} we find that 
$$\kk^{B_{b}^{\langle a-\varepsilon, b+\varepsilon \rangle}} \, \cong \, \kk^{B_{b}^{\langle a, b\rangle}}[\vec{\varepsilon}]. $$
Combining the last two isomorphisms with~\eqref{eq:Psiofcpctshift}, we find that 
\begin{equation}
 \Psi(\kk_{[a,b]}[-i] \star K_{\varepsilon}) \,\cong\, \Big(\kk^{B_{b}^{\langle a, b\rangle}}[-i]\Big)[\vec{\varepsilon}] \; \cong \;  \Psi\big(\kk_{[a,b]}[-i]\big)[\vec{\varepsilon}]
\end{equation}
using again Proposition~\ref{P:PsionIntervals} for the last isomorphism.
Similarly, in the case where $I$ is half-open, we obtain, for any $\varepsilon\geq 0$ and $i\in \mathbb{Z}$, 
\begin{equation}
 \Psi(\kk_{[a,b[} [-i] \star K_{\varepsilon}) \,\cong \, \Big(\kk^{B_{h}^{\langle a, b\rangle}}[-i]\Big)[\vec{\varepsilon}] \, \cong\, \Psi\big(\kk_{[a,b[} [-i]\big)[\vec{\varepsilon}] .
 \end{equation}
 \begin{equation}
 \Psi(\kk_{]a,b]} [-i] \star K_{\varepsilon})  \,\cong \, \Big(\kk^{B_{v}^{\langle a, b\rangle}}[-i]\Big)[\vec{\varepsilon}] \, \cong\,\Psi\big(\kk_{]a,b]} [-i]\big)[\vec{\varepsilon}].
\end{equation}

It remains to cover the case of an open interval $]a,b[$. Again by  Proposition~\ref{P:P38Berk18}, we have 
\begin{multline}\label{eq:PsiofcpctshiftBis}
 \Psi(\kk_{]a,b[}[-i] \star K_{\varepsilon}) \,\cong \,  \left\{ \begin{array}{ll}\Psi(\kk_{]a+\varepsilon,b-\varepsilon[}[-i]) & \mbox{ if }\varepsilon <\frac{b-a}{2} \\
      \Psi(\kk_{[b-\varepsilon,a+\varepsilon]}[-i]) & \mbox{ if }\varepsilon \leqslant\frac{b-a}{2} \end{array} \right .
      \\ \, \cong \, \left\{ \begin{array}{ll}S_i^{B_d^{\langle a+\varepsilon,b-\varepsilon\rangle}} & \mbox{ if }\varepsilon <\frac{b-a}{2} \\
      S_i^{B_b^{\langle b-\varepsilon,a+\varepsilon\rangle}} & \mbox{ if }\varepsilon \leqslant\frac{b-a}{2}                                                          \end{array}\right .
\end{multline} where the last isomorphism is given by proposition~\ref{P:PsionIntervals}.
Note that by definition the block $B_d$ is of type \textbf{dquad}$^+$ (see Lemma~\ref{L:Blockassociatedtointerval} and definition~\ref{D:blocksmodulesforMV}). Therefore as a persistent module over $\Delta_+$, we have, for $\varepsilon < \cfrac{b-a}{2}$, that
$$S_i^{B_{d}^{\langle a+\varepsilon, b -\varepsilon\rangle}} \,\cong\, \kk^{B_{d}^{\langle a-\varepsilon, b+\varepsilon \rangle}} [-i] \oplus \kk^{(B_{d}^{\langle a-\varepsilon, b+\varepsilon \rangle})^\dagger} [-i-1]$$ where the dual block $(B_{d}^{\langle a-\varepsilon, b+\varepsilon \rangle})^\dagger$ is of type \textbf{bb}$^{+}$.  
By Lemma~\ref{L:shiftofBlock} and Remark~\ref{R:caracofab} we find that 
$$\kk^{B_{d}^{\langle a-\varepsilon, b+\varepsilon \rangle}} [-i] \oplus \kk^{(B_{d}^{\langle a-\varepsilon, b+\varepsilon \rangle})^\dagger} [-i-1] \, \cong \,(\kk^{B_{d}^{\langle a, b \rangle}} [-i])[\vec{\varepsilon}] \oplus (\kk^{(B_{d}^{\langle a, b\rangle})^\dagger} [-i-1])[\vec{\varepsilon}]. $$
Combining these last two isomorphisms with~\eqref{eq:PsiofcpctshiftBis}, we find that,  for $\varepsilon < \cfrac{b-a}{2}$,
\begin{equation} \Psi(\kk_{]a,b[}[-i] \star K_{\varepsilon}) \,\cong \, S_i^{B_d^{\langle a, b\rangle}} [\vec{\varepsilon}]\, \cong \, \Psi(\kk_{]a,b[}[-i]) [\vec{\varepsilon}]
\end{equation}
as claimed.

\smallskip

It remains to consider the case  $\varepsilon \geq \cfrac{b-a}{2}$. We have still $\Psi(\kk_{]a,b[}[-i]) [\vec{\varepsilon}] \cong S_i^{B_{d}^{\langle a, b \rangle}} [\vec{\varepsilon}]$. As a graded persistent module over $\Delta_+$, by Lemma~\ref{L:shiftofBlock}, we have  
 that
$$S_i^{B_{d}^{\langle a, b \rangle}} [\vec{\varepsilon}] \,\cong\, \kk^{\big(B_{d}^{\langle a, b \rangle} -\vec{\varepsilon}\big)} [-i] \oplus \kk^{\big((B_{d}^{\langle a, b\rangle})^\dagger -\vec{\varepsilon}\big)} [-i-1].$$ 
But since  $\varepsilon \geq \cfrac{b-a}{2}$, we have that the death block  $\big(B_{d}^{\langle a, b \rangle} -\vec{\varepsilon}\big)$ is concentrated below the anti-diagonal $\Delta$, that is in $\mathbb{R}^2\setminus \Delta^+_{>0}$ and therefore $ \kk^{\big(B_{d}^{\langle a, b \rangle} -\vec{\varepsilon}\big)}\cong 0$.  Similarly, the birth block module  
$\big((B_{d}^{\langle a, b\rangle})^\dagger -\vec{\varepsilon}\big)$ is of type \textbf{bb}$^{-}$ precisely for $\varepsilon \geq \cfrac{b-a}{2}$. The infimum of the points included in this birth block has coordinates $(\varepsilon -b, a+\varepsilon)$. Therefore, 
$$\kk^{\big((B_{d}^{\langle a, b\rangle})^\dagger -\vec{\varepsilon}\big)} \, \cong \, \kk^{B_{b}^{\langle b-\varepsilon, a+\varepsilon \rangle}}.$$
Taking the (shifted) direct sum of this last two isomorphisms thus obtain, that, for $\varepsilon \geq \cfrac{b-a}{2}$, we have  
$$ S_i^{B_{d}^{\langle a, b \rangle}} [\vec{\varepsilon}] \,\cong\, 0 \oplus S_i^{B_{b}^{\langle b-\varepsilon, a+\varepsilon \rangle}}$$
and therefore the claim follows from the last case of~\eqref{eq:PsiofcpctshiftBis}.
\end{pf}

\begin{remark}
 Note that this kind of results has been proven with different assumptions in section 5 of \cite{BerkoukPetit}.
\end{remark}

\subsection{The isometry theorem between the interleaving distance on $\Delta^+$ and the graded bottleneck distance for sheaves} 
In this section, we will state and prove our main isometry theorem. 
Before that, we derive a few corollaries of the results we have obtained in sections~\ref{SS:Psi} and~\ref{SS:MVfunctor}.

\smallskip
Let $\mathcal{H}^*(-):  \D^b_{\R c}(\kk_\R) \to  \D^b_{\R c}(\kk_\R)$ be the endofunctor given by the cohohomology sheaf, that is, for any complex of sheaves  $F$,
by $$\mathcal{H}^*(F):=\bigoplus_{i\in \Z} H^i(F)[-i].
$$
\begin{cor}\label{C:MVcircPsiisId} Consider (the restrictions) $\Psi: \D^b_{\R c}(\kk_\R)\to \text{M-V}(\R)^{\s}$ and  
$\overline{(-)}^{MV}: \text{M-V}(\R)^{\s} \to  \D^b_{\R c}(\kk_\R)$. For any $F\in \D^b_{\R c}(\kk_\R)$, one has an isomorphism 
$$(\overline{~\cdot~ })^{\text{MV}} \circ \Psi (F)\, \cong \, F. $$
 Further, there is an  natural equivalence of functors 
 $(\overline{~\cdot~ })^{\text{MV}} \circ \Psi\; \simeq\; \mathcal{H}^*(-)$. 
\end{cor}
In other words, $\Psi$ is an natural section of the functor $\overline{(-)}^{MV}$ on strongly pointwise finite dimensional modules.  
\begin{pf} Let us prove the first claim.
 Since both functors $(\overline{~\cdot~ })^{\text{MV}}$ and $ \Psi $ commutes with shifts and direct sums 
 (propositions~\ref{P:PropertiesofMVSheafification} and \ref{P:PropertiesofPsi}), in view of the structure theorem~\ref{T:KSdecomposition},
 it is enough to construct the isomorphism for sheaves of the form $\kk_{I}$. 
 Now Proposition~\ref{P:PropertiesofMVSheafification}.2 (for $\varepsilon=0$) and Proposition~\ref{P:PsionIntervals} we have 
 \begin{equation}
\label{eq:MVcicrPSi0}   (\overline{\Psi(F) })^{\text{MV}} \circ \Psi (\kk_I) \;\cong \; \kk_I
 \end{equation}
 which is precisely giving such a claimed isomorphism  for an interval.
 
 \smallskip
 
 Let us prove the natural equivalence; we will denote by $(-)^*$ the linear dual as before. 
 By definition of $\Psi(F)$, if  $F\in \D^b_{\R c}(\kk_\R)$, one has, 
\begin{equation}\label{eq:MVcicrPSi1}(\overline{\Psi(F) })^{\text{MV}} \circ \Psi (F) \, = \, 
 \ell \left(\bigoplus_{i\in \Z} \overline{\Big(H^i\big(\Gamma(I^\bullet_{|\, (-)})\big)\Big)^{*}}[-i]\right)
 = \ell \left(\bigoplus_{i\in \Z} \widetilde{\Big(H^i\big(\Gamma(I^\bullet_{|\, (-)})\big)\Big)^{*}}[-i]\right). \end{equation}
 By definition~\ref{D:defofBar} and lemma~\ref{L:valueoftilde}, denoting $(-)^a$ the sheafification functor, we have, for any open set $U\subset \R$, 
 \begin{equation}\label{eq:MVcicrPSi2}
  \overline{\Big(H^i\big(\Gamma(I^{\bullet}_{|\, U})\big)\Big)^{*}} = \left(\widetilde{\Big(H^i\big(\Gamma(I^{\bullet}_{|\, U})\big)\Big)^{*}}\right)^{a}
  =\left(\varprojlim_{]-x,y[\subset U} {\Big(H^i\big(\Gamma(I^\bullet_{|\, ]-x,y[})\big)\Big)^{**}}\right)^{a}.
 \end{equation} Recal that $\Gamma (F_{|V})=F(V)$ by definition and therefore
 the restriction homomorphisms $I^\bullet(U) \to I^\bullet (]-x, y[)$ yields the canonical morphism 
 \begin{equation}\label{eq:MVcicrPSi3}
  \varphi_U: H^i(I^\bullet(U)) \hookrightarrow\Big(H^i\big(I^\bullet (U)\big)\Big)^{**} \cong \Big(H^i\big(\Gamma(I^\bullet_{|\, U})\big)\Big)^{**}
  \to \hspace{-0.5pc}\varprojlim_{]-x,y[\subset U} \hspace{-0.5pc}{\Big(H^i\big(\Gamma(I^\bullet_{|\, ]-x,y[})\big)\Big)^{**}}
 \end{equation}
where the left map is the canonical homorphism for a vector space to its bidual. Since the cohohomology $H^i(F)$ of a a complex of sheaf is
the sheafification of 
the presheaf $U\mapsto H^i(F(U))$, combining~\eqref{eq:MVcicrPSi1},~\eqref{eq:MVcicrPSi2} and \eqref{eq:MVcicrPSi3} we obtain a morphism of sheaves:
\begin{multline}
 \mathcal{H}^*(F)(U) =  \bigoplus_{i\in \Z} H^i(F)(U)[-i]  \; \cong \;  
 \bigoplus_{i\in \Z} H^i(I^\bullet)(U)[-i] \\ \xymatrix{ \ar[rr]^{ \bigoplus (\varphi_U)^{a}\qquad\quad} &&
 (\overline{\Psi(F) })^{\text{MV}} \circ \Psi (F)(U)}
\end{multline}
where the first isomorphism is given by the quasi-isomorphism $F\stackrel{\simeq}\to I^\bullet$. To see that this map is an isomorphism, it is enough to check it for the stalks, 
and therefore, since
  that map commutes with finite direct sum and shifts (a finite direct sum of injective resolutions is an injective resolution), 
  it is enough to prove that this map is an isomorphism when $F=\kk_I$. The latter follows from~\eqref{eq:MVcicrPSi0}.
\end{pf}
\begin{cor}\label{C:fullorfiathful}
The functor  $\overline{(-)}^{MV}: \text{M-V}(\R)^{\s} \to  \D^b_{\R c}(\kk_\R)$ is essentially surjective. 
\end{cor}
\begin{pf}
For  any    constructible sheaf $F^\bullet$, $\Psi(F^\bullet)$ is a 
strongly pointwise finite dimensional Mayer-Vietoris system and Corollary~\ref{C:MVcircPsiisId} gives a natural isomorphism 
$F^\bullet \cong \overline{(\Psi(F^\bullet))}^{MV} $. Therefore, $F^\bullet$ is in the essential image of $\overline{(-)}^{MV}$.  
\end{pf}

\begin{cor}\label{C:MVdisO}
 Let $M\in \mbox{MV}(\R)^\s$ be a strongly pointwise finite dimensional Mayer-Vietoris system. Then 
 $$d_I^{\text{MV}}(M, \Psi ( \overline{M}^{MV})) = 0. $$
\end{cor}
In other words, though $\Psi \circ \overline{(-)}^{MV}$ is not an equivalence, it maps an object to an object which is at distance $0$ fom itself.
\begin{pf}
 By statement 3 of Proposition~\ref{P:PropertiesofMVSheafification}, it is sufficient to prove that $\overline{M}^{MV}$ and $\overline{\Psi(\overline{M}^{MV})}^{MV}$ are isomorphic in $\D^b_{\R c}(\kk_\R)$. But Corollary~\ref{C:MVcircPsiisId} implies  $\overline{\Psi(\overline{M}^{MV})}^{MV} \cong \overline{M}^{MV}$ and the result follows.
\end{pf}

We can now state our isometry theorem
\begin{thm}[Isometry]\label{T:MainIsometry} The Mayer-Vietoris sheafification functor and the functor $\Psi$ are pseudo-isometries between the interleaving distance and the convolution and bottleneck distances for sheaves. That is, for all $M, N \in \mbox{MV}(\R)^\s$,  one has equalities 
$$d_I(M, N) = d_C (\overline{M}^{MV}, \overline{N}^{MV}) = d_B(\mathcal{B}(\overline{M}^{MV}), \mathcal{B}(\overline{N}^{MV})). $$ 
And for all constructible sheaves $F, G \in \D^b_{\R c}(\kk_\R)$, one has 
$$ d_B(\mathcal{B}(F), \mathcal{B}(G)) = d_C(F,G) =d_I(\Psi(F), \Psi(G)).$$
\end{thm}
\begin{pf}Theorem~\ref{T:DerivedIsometry} implies already the equality between bottleneck and convolution distances. 

 By Proposition~\ref{P:PropertiesofMVSheafification}, the Mayer-Vietoris sheafification functor $\overline{(-)}^{MV}:  
 \mbox{MV}(\R)^\s\to \D^b_{\R c}(\kk_\R)$ maps the shift functor $[\vec{\varepsilon}]$ onto the convolution functor $(-)\star K_{\varepsilon}$  
 and therefore if  $M, N \in \mbox{MV}(\R)^\s$ are $\varepsilon$-interleaved,
 then $\overline{M}^{MV} \sim_{\varepsilon} \overline{N}^{MV}$ in $ \D^b_{\R c}(\kk_\R)$. Thus,  for all $M, N \in \mbox{MV}(\R)^\s$,  one has 
 \begin{equation}
  \label{eq:inegMVtoShv}  d_C (\overline{M}^{MV}, \overline{N}^{MV}) \, \leqslant \, d_I(M, N).
 \end{equation}
Similarly,  Proposition~\ref{P:Psipreservesinterl} implies that $\Psi$ sends the convolution  $(-)\star K_{\varepsilon}$ functor to the shift functor and thus, we also have that, for all $F,\, G\in  \D^b_{\R c}(\kk_\R)$, one has  
\begin{equation}
  \label{eq:inegSHtoMV} d_I(\Psi(F), \Psi(G)) \, \leqslant \, d_C (F, G). 
 \end{equation}
 From~\eqref{eq:inegMVtoShv} and~\eqref{eq:inegSHtoMV} we get, for all $M, N \in \mbox{MV}(\R)^\s$, that
 \begin{equation}\label{ineq:chained}
  d_I(\Psi\big(\overline{M}^{MV} \big), \Psi\big(\overline{N}^{MV} \big)) \, \leqslant \,
  d_C (\overline{M}^{MV}, \overline{N}^{MV})
  \, \leqslant \, d_I(M, N).
 \end{equation}
The triangular inequality and Corollary~\ref{C:MVdisO} implies
\begin{multline}\label{ineq:Icomposed}
 d_I(M,N) \leqslant  d_I(\Psi\big(M, \Psi\big(\overline{M}^{MV} \big))
+ d_I(\Psi\big(\overline{M}^{MV} \big), \Psi\big(\overline{N}^{MV} \big)) 
+ d_I(\Psi\big(\overline{N}^{MV} \big), N)\\ 
= d_I(\Psi\big(\overline{M}^{MV} \big), \Psi\big(\overline{N}^{MV} \big))
 \end{multline}
Combining inequality~\eqref{ineq:Icomposed} with~\eqref{ineq:chained}, we ontain that all inequalities in~\eqref{ineq:chained} are egalites which gives the first claim 
$$d_I(M, N) = d_C (\overline{M}^{MV}, \overline{N}^{MV}).$$
To prove the remaining one, we use Corollary~\ref{C:MVcircPsiisId}. 
This gives us, for any $F, \, G\in  \D^b_{\R c}(\kk_\R)$ isomorphisms 
$F\cong \overline{\Psi(F)}^{MV}$ and $G\cong \overline{\Psi(G)}^{MV}$ and therefore we have
\begin{equation}\label{eq:lastegality}
 d_C(F,G) \, = \, d_C\big( \overline{\Psi(F)}^{MV}, \overline{\Psi(G)}^{MV}\big) \, =
 d_I(\Psi(F), \Psi(G))
\end{equation}
since we just proved that $\overline{(-)}^{MV}$ is an isometry. The equality~\eqref{eq:lastegality} concludes the proof of the theorem.
\end{pf}

In particular the theorem allows to compute the bottleneck or convolution distance for sheaves using interleaving distance for persistence modules and vice-versa. Furthermore, we recover as a corollary the following result of \cite{KS18}.
\begin{cor} If $X$ is a locally contractible compact topological manifold, and $f,g : X \to \R$ are continuous constructible functions,  one has : $$d_C(\Rr f_* \kk_X, \Rr g_* \kk_X) \leq \sup_{x\in X} \|f(x) - g(x) \|$$
\end{cor}
\begin{pf}
 Theorem~\ref{T:MainIsometry} and proposition~\ref{P:MVmapstoRf} imply that $$d_C(\Rr f_* \kk_X, \Rr g_* \kk_X)= d_I(M^f, M^g)\leqslant  \sup_{x\in X} \|f(x) - g(x) \|$$ 
 by Proposition~\ref{P:StabilityfordMV}.
\end{pf}

\subsection{A detailled example} \label{SS:Projectionsfromcircle}
The following example was suggested to us by Justin Curry. 

Let $\mathbb{S}^1 = \{(x,y) \in \R^2 \mid x^2 + y^2 = 1 \}$ be the  circle embedded in $\R^2$. 
Let  $f : \mathbb{S}^1 \to \R$ be the first coordinate projection and $p : \mathbb{S}^1 \to \R $ 
be the constant map with value zero.
\begin{figure}
\begin{center}
\includegraphics[scale=0.8]{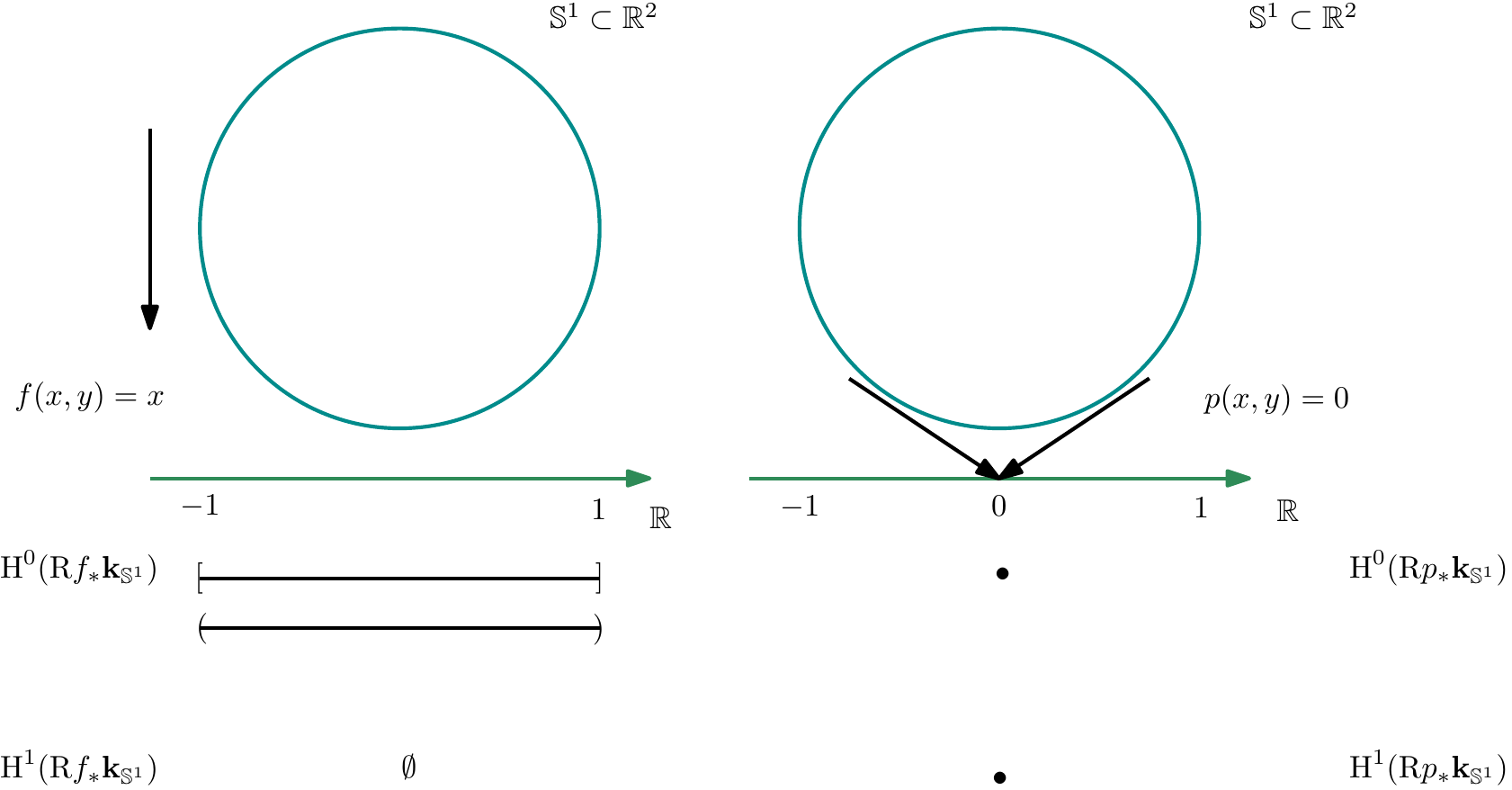}
\label{fig:circleprojection}
\caption{The map $f$ and $p$ together with the intervals on which are supported the degree $0$ and $1$ part of the associated sheaves, see~\eqref{eq:RfandRpofcircle} and proposition~\ref{P:exofcircleproj}.}
\end{center}
\end{figure}
From example~\ref{Ex:MVfromFct} we obtain two Mayer-Vietoris systems $M^f$ and $M^p$, which are given, for any $(x,y)\in \Delta^+$, by
$M^f(x,y)= H_*(f^{-1}(]-x,y[))$ and $M^p(x,y)= H_*(p^{-1}(]-x,y[))$. 

Using the same notation as in Lemma~\ref{L:Blockassociatedtointerval} we have.
\begin{prop}\label{P:exofcircleproj}
 One has $$M^f \, \cong\, S_0^{B_b^{]-1,1[}}\oplus S_0^{B_d^{[-1,1]}}, \qquad 
 M^p \, \cong\, S_0^{B_b^{]0,0[}}\oplus S_1^{B_b^{]0,0[}}.$$
 
 In particular,  $\overline{M^f}^{MV} =\kk_{(-1,1)} \oplus \kk_{[-1,1]}$ and  $\overline{M^p}^{MV} = \kk_{\{0\}} \oplus \kk_{\{0\}}[-1]$. Furthermore, $\Psi(\kk_{(-1,1)} \oplus \kk_{[-1,1]}) \cong M^f$ and $\Psi (\kk_{\{0\}} \oplus \kk_{\{0\}}[-1]) \cong M^p$.
\end{prop}
\begin{figure}
\begin{center}
\includegraphics[scale=0.8]{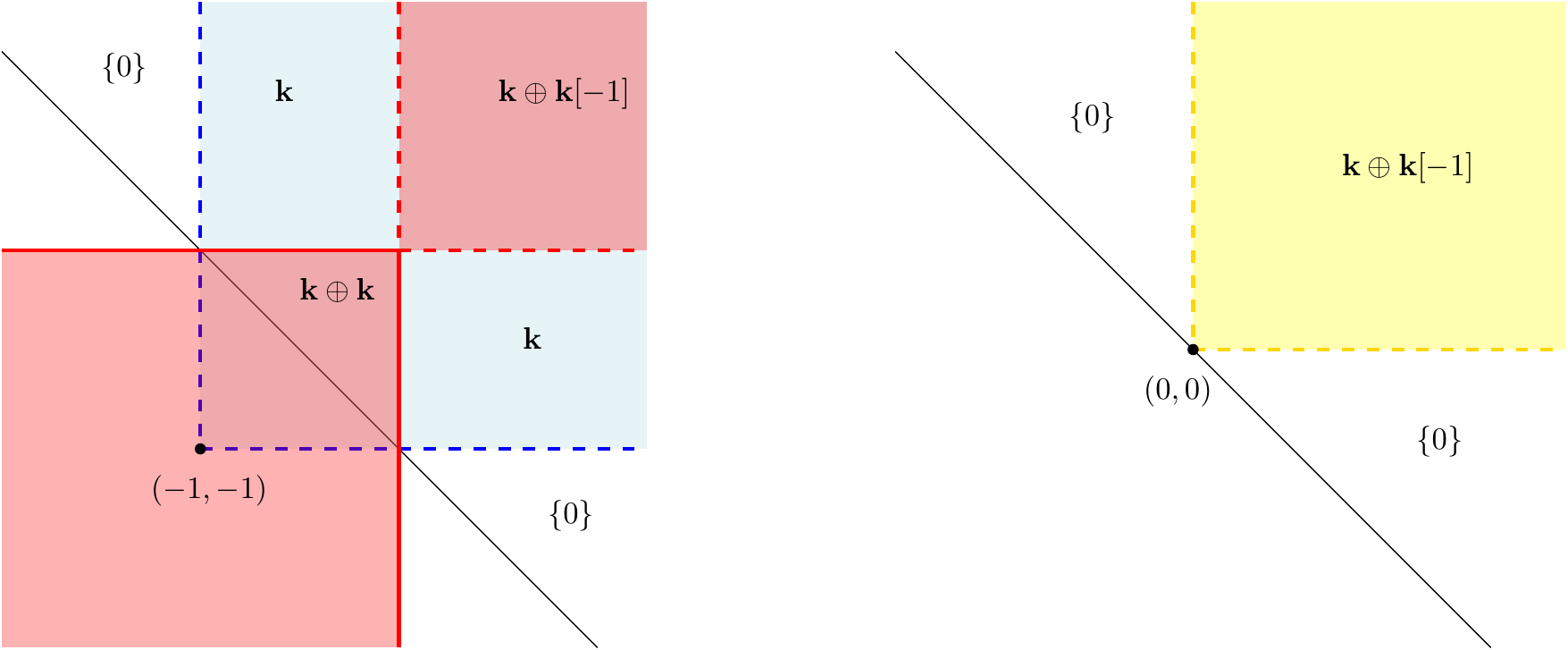}
\label{fig:moduleforcircle}
\caption{On the left, the value of the  MV system $M^f$ where the blue part is a birth block and the red part are a death block and its dual. On the right, the value of the MV system $M^p$ where the yellow part refers to the (reunion of) two birthblocks. 
The dashed lines pictures boundary which are not inside the blocks. }
\end{center}
\end{figure}
\begin{pf}
 The preimage of $f$ satisfies 
 $$f^{-1}(]-x,y[) =\left\{ \begin{array}{ll} \emptyset & \mbox{if } -x\geqslant 1 \mbox{ or } y\leqslant -1\\
 S^1 & \mbox{if } -x<-1 \mbox{ and } y>1 \\
 \mbox{two intervals} & \mbox{if } -1\leqslant -x <y \leqslant 1\\
 \mbox{one interval} & \mbox{if } -x<-1<y\leqslant 1 \mbox{ or } -1\leqslant -x <1<
 y.
 \end{array}\right. $$
 This gives that $M^f$ has the module decomposition given in figure~\eqref{fig:moduleforcircle} which is exactly the decomposition of $M^f$ into a \textbf{bb}$^{-}$ module with infimum $(-1,-1)$ in degree $0$ and a module associated to the deathblock with supremum $(1,1)$. 
 The image of $\Psi$ is  given by adidtivity and Proposition~\ref{P:PsionIntervals}:
 \begin{multline*}\Psi(\Rr f_*\kk_{\mathbb{S}^1} ) \cong \Psi( \kk_{(-1,1)} \oplus \kk_{[-1,1]}) \cong \psi(\kk_{(-1,1)}) \oplus \Psi(\kk_{[-1,1]}) \cong 
 S_0^{B_b^{]0,0[}}\oplus S_1^{B_b^{]0,0[}} \cong M^f.
 \end{multline*}
 Applying Corollary~\ref{C:MVcircPsiisId} (or using Proposition~\ref{P:PropertiesofMVSheafification} directly) yields $\overline{M^f}^{MV} =\Rr f_* S^1$.
Similarly, we have that 
$$p^{-1}(]-x,y[) =\left\{ \begin{array}{ll} S^1 & \mbox{if } -x<0<y \\
                           \emptyset & \mbox{else}.
                          \end{array}\right.
 $$
 and thus $M^p \, \cong\, S_0^{B_b^{]0,0[}}\oplus S_1^{B_b^{]0,0[}}$ as can be seen in figure~\eqref{fig:moduleforcircle} as well.
We apply again Proposition~\ref{P:PsionIntervals} and Corollary~\ref{C:MVcircPsiisId} to conclude.
\end{pf}
In particular we recover the computation of~\cite{Berk18} for the derived images sheaves 
$\Rr f_*\kk_{\mathbb{S}^1}$  and $\Rr p_*\kk_{\mathbb{S}^1}$ :
\begin{eqnarray}\label{eq:RfandRpofcircle}
 \Rr f_*\kk_{\mathbb{S}^1} \,\cong\, \kk_{]-1,1[} \oplus \kk_{[-1,1]}, \qquad 
 \Rr p_* \kk_{\mathbb{S}^1} \,\cong \,\kk_{\{0\}} \oplus \kk_{\{0\}}[-1].
\end{eqnarray}
Furthermore,  we can find $1$-interleaving between $S_0^{B_d^{[-1,1]}} $ and $S_1^{B_b^{]0,0[}}$ as well as $1$-inteleaving between $S_0^{B_b^{]-1,1[}}$ and $S_0^{B_b^{]0,0[}}$. Therefore,
the decomposition of the proposition gives a $1$-inteleaving for $M^f$ and $M^p$.

\section*{Acknowledgments}
The authors wish to thank Benedikt Fluhr for his very useful comments, notably 
the non-naturality of the isomorphism of Lemma~\ref{L:DefadnFunctorialityPsi}, and for providing the counter-example of Remark~\ref{rk:Psinotfaithful}.

\bibliographystyle{alpha}

\end{document}